\documentclass[11pt]{article}
\usepackage[a4paper,margin=2.5cm]{geometry}
\usepackage{amsmath,amssymb,amsthm,mathtools}
\usepackage{booktabs,array,graphicx,enumitem,xcolor}
\usepackage[hidelinks]{hyperref}
\usepackage[capitalise,nameinlink]{cleveref}
\newtheorem{theorem}{Theorem}[section]
\newtheorem{proposition}[theorem]{Proposition}
\newtheorem{lemma}[theorem]{Lemma}
\newtheorem{corollary}[theorem]{Corollary}
\theoremstyle{definition}
\newtheorem{definition}[theorem]{Definition}
\newtheorem{assumption}[theorem]{Assumption}
\newtheorem{remark}[theorem]{Remark}
\newtheorem{conjecture}[theorem]{Conjecture}
\numberwithin{equation}{section}

\newcommand{\R}{\mathbb{R}}
\newcommand{\Z}{\mathbb{Z}}
\newcommand{\E}{\mathbb{E}}

\newcommand{\D}{\mathcal{D}}
\newcommand{\M}{\mathcal{M}}
\newcommand{\norm}[1]{\lVert #1\rVert}
\newcommand{\abs}[1]{\lvert #1\rvert}
\newcommand{\ip}[2]{\langle #1,#2\rangle}
\newcommand{\eps}{\varepsilon}
\newcommand{\hatt}[1]{\widehat{#1}}
\newcommand{\sgn}{\operatorname{sgn}}
\newcommand{\Id}{\mathrm{Id}}
\newcommand{\BMO}{\mathrm{BMO}}
\newcommand{\curv}{\nu_\tau}
\DeclareMathOperator{\supp}{supp}

\title{\textbf{Random Dyadic Shift Scattering:}\\[3pt]\Large\textbf{Scale-Uniform Deformation Stability\\ and Phase-Aligned Hilbert Multipliers}}
\author{Daewon Chung\\[2pt]\small Department of Mathematics, Keimyung University, Daegu, Korea\\\small\texttt{dwchung@kmu.ac.kr}}
\date{\small Preprint, November 2026}

\begin{document}
\maketitle

\begin{abstract}
\noindent
Petermichl's theorem represents the Hilbert transform by averaging a single dyadic shift over dilations and random dyadic lattices. We study how the deformation stability of the lattice-averaged shift, and of a parameter-free feature map $\Phi$ built from it in the spirit of Mallat's wavelet scattering, depends on the number $J+1$ of scales. For the Haar generator as well as compactly supported $C^1$ orthonormal wavelets, the covariance error under a deformation $f\mapsto f(\cdot-\tau(\cdot))$ is $O((J+1)\norm{\tau'}_\infty\norm f)$ without a $\norm\tau_\infty$ term, because lattice averaging absorbs the translation part of $\tau$ exactly. The deformation defect splits into a first-order Calder\'on commutator, bounded uniformly in $J$ by the $T(1)$ theorem, and a second-order remainder. On a single lattice, affine warps and $C^{1,\delta}$ generators give a lower bound $c\,\sigma_\eps\sqrt{J+1}-C\eps$ for the defect, with an explicit constant $\sigma_\eps\ge0$ computed to be positive for sym4 and db4 at three values of $\eps$. Averaging over lattices turns the defect into a Fourier-multiplier quadratic form, yielding $J$-uniform bounds both for affine warps, with a logarithmic loss in the Haar case, and for $C^2$ warps with $C^{1,\delta}$ generators, and hence a scale-uniform Mallat-type stability estimate for $\Phi$. With the standard child pair, the Hilbert component of the averaged shift nearly cancels for smooth Daubechies-type generators; re-indexing the children by integers restores this component. For Haar and general warps, numerical computations agree with a conjectured $J$-uniform bound.
\end{abstract}

\noindent\textbf{Keywords:} dyadic shift; Hilbert transform; random dyadic lattices; deformation stability; Calder\'on commutator; scattering transform\\
\textbf{MSC 2020:} 42B20, 42C40, 42A45 (primary); 68T07 (secondary)

\tableofcontents
\newpage

%==================================================================
\section{Introduction}\label{sec:intro}

\paragraph{Dyadic shifts.}
Petermichl [P00] established that averaging a single \emph{dyadic shift} over random dyadic lattices and their dilations yields, up to a constant multiple, the Hilbert transform $H$. This dyadic shift  maps the Haar coefficient of each interval to the difference of its two children. The extension to every Calder\'on--Zygmund operator came through Hyt\"onen's dyadic representation theorem [H12]; subsequently, Hyt\"onen--Lappas [HL22] and Di~Plinio--Wick--Williams [DWW23] substituted smooth compactly supported wavelets for the Haar system. Each of these results decomposes a singular integral into simple, local, multiscale pieces, with averaging over random lattices restoring translation and dilation symmetry.

\paragraph{The network.}
We use this decomposition to construct a network. For the \emph{dyadic shift scattering network} (DSSN), a first-order channel takes the form
\[
\Psi_1f=A_J\,\Psi_1'f,\qquad\Psi_1'f=\E_\beta\abs{S^J_\beta f},
\]
where the dyadic shift $S^J_\beta$ acts on the random lattice $\D^\beta$ and is truncated to $J+1$ scales; $\abs{\cdot}$ denotes the pointwise modulus, while $A_J$ denotes a low-pass average at scale $2^J$ (\cref{sec:setting}). As in Mallat's wavelet scattering [M12, BM13], the construction combines a linear multiscale operator with a modulus and a low-pass filter. Two distinctions remain: instead of a convolution, the linear stage uses a lattice-dependent \emph{shift}; instead of obtaining invariance from the filter bank, it obtains it by averaging over lattices. With no trained parameters, this network serves as our concrete model for four questions.
\begin{enumerate}[label=(Q\arabic*),leftmargin=2.5em,nosep]
\item Which Fourier multiplier is realized by the lattice-averaged shift, and how is it affected by whether averaging includes translations alone or dilations as well?
\item Under a deformation $L_\tau f(x)=f(x-\tau(x))$, what decomposition separates the defect of the shift into a first-order commutator and higher-order remainders, and how does each depend on the number $J+1$ of scales?
\item Why is accumulation over scales controlled by averaging over lattices, and is one lattice sufficient?
\item How does positioning the children of a smooth generator influence the phase of the averaged shift and, consequently, its Hilbert component?
\end{enumerate}
Questions (Q2) and (Q3) provide an analogue of Mallat's deformation stability theorem [M12], with refinements due to Wiatowski--B\"olcskei [WB18] and, under minimal regularity of $\tau$, Nicola--Trapasso [NT23]. Convolutional filter banks are the subject of those works; the lattice-based network instead relies on a different mechanism. If the hypotheses in \cref{mt:B}(ii)--(iii) hold, the stability constant for the \emph{averaged} network is uniform in $J$. For a single lattice, however, this uniformity fails whenever the constant $\sigma_\eps$ in \cref{thm:single-lower} is positive.

\paragraph{What is known.}
Question (Q1) already has an essential answer in the literature: averaging the shift over translations and dilations gives a multiple of $H$ [P00, H08]. Our \cref{prop:mult,thm:haarkappa} express this fact on the Fourier side through the explicit symbol $m_0$ used throughout the remainder of the paper. Because translation averaging by itself commutes only with dyadic dilations, the associated symbol is log-periodic. The resulting ripple is quantified by \cref{cor:transfail} (about $\pm10\%$ per octave for Haar). This is practically relevant, since the network must choose the number of dilations to sample. Questions (Q2)--(Q4) are the subject of the new results.

\paragraph{Main results.}
Write $\eps=\norm{\tau'}_\infty\le\frac12$ and $\curv=\norm{\tau''}_\infty$, using the coarsest cell as the unit of measurement.

\begin{theorem}[Finite-$J$ covariance]\label{mt:A}
Let $\psi$ denote either the Haar function or a compactly supported $C^1$ orthonormal wavelet, allowing any child offset $m\in\{0,\dots,L-1\}$, and take $\eps\le\frac12$. Then
\[
\norm{\Psi_1'(L_\tau f)-L_\tau\Psi_1'f}_2\le C_\psi(J+1)\,\eps\,\norm f_2 ,
\]
with $C_h\le10.4$ in the Haar case ($C_h\to\sqrt{48}\approx6.9$ as $J\to\infty$). The quantity $\norm\tau_\infty$ does not enter. For a single Haar lattice, the corresponding defect has exact order $\eps^{1/2}$ (\cref{prop:single-haar}).
\end{theorem}
\noindent The proof of \Cref{mt:A} appears in \cref{sec:finiteJ} as \cref{thm:cov}.

\begin{theorem}[Single lattice versus lattice average]\label{mt:B}
Let $G_\beta$ denote the Jacobian-free deformation defect for the shift on lattice $\D^\beta$ (\cref{sec:finiteJ}).
\begin{enumerate}[label=(\roman*),nosep]
\item For every lattice and for $\psi\in C^{1,\delta}$, we have $\norm{G_\beta}\le C(\eps+(J+1)\eps^{1+\delta})$. A bound of $C\eps$ holds for the first-order part, uniformly in $J$. In the case of affine warps, the remainder obeys the lower bound $\norm{G_\beta}\ge c\,\sigma_\eps\sqrt{J+1}-C\eps$ for all $J\ge J_0(\eps)$ and every lattice; here $\sigma_\eps\ge0$ denotes an explicit lacunary-series constant (\cref{thm:single-lower}). Consequently, a single lattice is not $J$-uniform \emph{whenever $\sigma_\eps>0$}. Our computations give $\sigma_\eps>0$ for aligned sym4 and db4 at $\eps\in\{0.05,0.1,0.2\}$; a proof of $\sigma_\eps>0$ for every admissible generator and every $\eps$ is not provided.
\item For affine warps, the bounds are $\sup_f\E_\beta\norm{G_\beta f}^2/\norm f^2\le C\eps^2$ for $C^{1,\delta}$ generators and $\le C\eps^2\log^2(e/\eps)$ for Haar, with uniformity in $J$.
\item For $C^2$ warps with $\psi\in C^{1,\delta}$, the bound $(\E_\beta\norm{G_\beta f}^2)^{1/2}\le C(\eps+\eps\curv^\delta+\eps^\delta\curv)\norm f$ holds uniformly in $J$.
\end{enumerate}
\end{theorem}
\noindent In \cref{mt:B}, \cref{thm:A,thm:single-lower} are combined in part (i), whereas parts (ii) and (iii) are \cref{thm:B,thm:C}.

\begin{corollary}[Scale-uniform stability]\label{mt:C}
For $\psi\in C^{1,\delta}$ and $\tau\in C^2$,
\[
\norm{\Phi(L_\tau f)-\Phi f}\le C_\psi\Big[\eps+\eps\curv^\delta+\eps^\delta\curv+2^{-J}\norm{\phi'}_1\norm\tau_\infty\Big]\norm f_2,
\]
with $C_\psi$ not depending on $J$.
\end{corollary}
\noindent \Cref{mt:C} is given by \cref{thm:main}.

\paragraph{Smooth generators and phase alignment.} Although the scale-uniform theory requires $C^{1,\delta}$ generators, the placement of the two children relative to their parent determines the Hilbert component of the averaged shift. With the standard child pair $(2k,2k+1)$ for a Daubechies-type wavelet, almost complete cancellation occurs in this component: compared with $1.77$ for Haar, the Monte Carlo constant $\abs{c_0}=1/\abs{\operatorname{Im}\kappa_\psi}$ equals $43$ for db4 and $1851$ for coif3. The phase of the parent--child cross-spectrum accounts for this cancellation (\cref{ssec:smooth}). When the child offset is selected to maximize $\abs{\operatorname{Im}\kappa_\psi}$, $\abs{c_0}$ returns to the range between $1.90$ and $2.56$ for db2, db4, sym4, sym6 and coif3 (\cref{ssec:smooth}, \cref{tab:aligned}). Every value is calculated using the exact spectrum rather than sampling. Under their regularity hypotheses on $\psi$, \Cref{mt:A,mt:B,mt:C} hold for any child offset $m\in\{0,\dots,L-1\}$; alignment may therefore be selected solely for the Hilbert component. Since db2 is not $C^1$, only the multiplier results apply to it.

\paragraph{Extensions.} Three practically used network variants also admit these results (\cref{ssec:ext}). (1) With \emph{scale-separated} channels $\E_\beta\abs{S_{\beta,j}f}$, one at each scale as in scattering, $\sqrt{J+1}$ replaces the finite-depth factor $J+1$ as an improvement. The affine and $C^2$ bounds continue to be $J$-uniform, while no $T(1)$ theorem is needed for the first-order part (\cref{thm:sep}). (2) Non-expansiveness holds for a \emph{second layer}, which doubles the stability constant (\cref{cor:depth2}). (3) When \emph{dilation averaging} is added, as required by the Hilbert-transform identity, $\eps$ stays unchanged, while $\curv$ increases by at most a factor of two (\cref{cor:dil}). For each of these three variants, numerical experiments appear in \cref{ssec:num-ext}.

\noindent For Haar with non-affine warps, our conjectured bound is $C(\eps\log\frac e\eps+\curv)$ (\cref{conj:haar}). This is consistent with numerical computations: when $J\le8$ and $\eps\in[0.02,0.3]$, the smooth warp's averaged worst-case constant $\sigma_{\max}/\eps$ saturates as a function of $J$ at approximately $2.0$--$2.2$ (\cref{ssec:num-worst}). For $J\le7$, the reported values come from Lanczos and match the largest eigenvalue of the full matrix wherever that matrix was assembled; for $J=8$, they are lower bounds obtained by power iteration. Neither provides an upper bound for the operator norm.

\paragraph{The mechanism.}
Three facts underlie the proofs. (a) \emph{Local translation identity} (\cref{lem:local}): by converting the deformed output at $x$ to the output associated with the input translated by $\tau(x)$, the lattice average absorbs the translation part of $\tau$ exactly, leaving only the local distortion. (b) \emph{Commutator structure}: a change of variables expresses the residual defect through the kernel $\sum_j[k_j(x,v)-k_j(x,v+\tau(x)-\tau(v))]$, which is a dyadic relative of Calder\'on's first commutator [Cal65]. Uniformly in $J$, its first-order Taylor part meets the hypotheses of the $T(1)$ theorem. (c) \emph{Averaging turns lattice operators into multipliers} (\cref{lem:avg}): averaging $\norm{T_\beta g}^2$ with respect to $\beta$ yields a Fourier multiplier quadratic form. Away from $\abs\xi\approx2^j$, the scale-$j$ symbol is small; consequently, only boundedly many scales contribute at any given frequency, allowing frequency-by-frequency control of the sum over scales. By contrast, a single-lattice operator is not a Fourier multiplier; in that case, \cref{thm:single-lower} demonstrates that the scale sum can grow.

\paragraph{Related work.}
\emph{Dyadic shifts and representation theorems.} Random dyadic lattices originate in the work of Nazarov--Treil--Volberg [NTV03]. In Petermichl's representation of $H$ [P00], the average is taken over translations and dilations. The same approach yielded the Riesz transforms for Petermichl--Treil--Volberg [PTV02], every sufficiently smooth one-dimensional convolution kernel for Vagharshakyan [Vag10], and every Calder\'on--Zygmund operator for Hyt\"onen [H12]. Strong convergence of Petermichl's average was proved by Hyt\"onen [H08], who also established, for the continuous family of shifts, the identity $\int_0^\infty k_t\,dt/t=-\frac\pi8H$ equivalent to our Haar constant $\kappa_h$. In the representation theorems of Hyt\"onen--Lappas [HL22] and Di~Plinio--Wick--Williams [DWW23], smooth wavelets appear via averaged shifts of all complexities. To our knowledge, neither the Hilbert component of the averaged complexity-one shift for smooth generators nor its dependence on the child offset has been treated; nor has deformation stability for any of these lattice operators. A different complexity-one shift, the \emph{dyadic Hilbert transform}, was introduced by Domelevo--Petermichl [DP22, DP26] and rotates each interval's two children. Rather than averaging, the comparison with $H$ uses $L^p$ norm inequalities for Banach-valued functions; the lattice average of this shift is a zero multiple of $H$ (\cref{rem:sib}).

\emph{Commutators.} The deformation defect is related dyadically to Calder\'on's first commutator [Cal65]. The classical theory of that commutator encompasses Coifman--Meyer [CM75] as well as the $T(1)$ theorem of David--Journ\'e [DJ84]; we apply the latter to the defect's first-order part.

\emph{Scattering.} Wavelet scattering, together with its deformation stability, is developed in [M12, BM13, WB18, NT23]. On a fixed multiresolution, Haar scattering on graphs [CCM14] employs the same additions, subtractions and moduli. The procedure of averaging over lattice offsets is cycle spinning [CD95]. Our alignment condition recalls the half-sample delay condition for wavelet Hilbert pairs [S01] underlying the dual-tree complex wavelet transform [SBK05]. Whereas two filter banks produce the Hilbert pair there, lattice averaging of one real wavelet gives rise to a Hilbert component here.

\paragraph{Organization.}
\Cref{sec:setting} establishes the notation. The averaged multiplier is computed in \Cref{sec:mult}, which also treats aligned smooth generators. The finite-$J$ covariance bounds are proved in \Cref{sec:finiteJ}. Both the single-lattice and lattice-averaged analysis, together with the scale-uniform stability theorem, appear in \Cref{sec:uniform}. Numerical experiments are presented in \Cref{sec:num}; limitations and further directions are discussed in \cref{sec:disc}. The code and data for reproducing the computations are listed in \cref{app:repro}.

%==================================================================
\section{Setting}\label{sec:setting}

We adopt $\hatt f(\xi)=\int_\R f(x)e^{-2\pi ix\xi}\,dx$; consequently, $\widehat{Hf}(\xi)=-i\sgn(\xi)\hatt f(\xi)$ and $T_af(x)=f(x-a)$. Unless otherwise specified, every norm is an $L^2(\R)$ norm.

\subsection{Generators, child offsets and shifts}
\begin{definition}[Generator, child offset, shift]\label{def:shift}
Let $\psi$ denote a real orthonormal wavelet satisfying $\supp\psi\subset[0,L]$ for some integer $L\ge1$. For the Haar function $h=\mathbf 1_{[0,1/2)}-\mathbf 1_{[1/2,1)}$, one has $L=1$. Choose a fixed integer \emph{child offset} $m\in\{0,\dots,L-1\}$ and define
\[
\Theta(s)=\Theta^{(m)}(s)=\psi(2s-m)-\psi(2s-m-1),\qquad\supp\Theta\subset\big[\tfrac m2,\tfrac{m+L+1}2\big]\subset[0,L].
\]
For an interval $I=2^{-j}([0,1)+k)+\gamma$, set $\psi_I(x)=2^{j/2}\psi(2^j(x-\gamma)-k)$ and $\Theta_I(x)=2^{j/2}\Theta(2^j(x-\gamma)-k)$; equivalently, $\Theta_I=\frac1{\sqrt2}(\psi_{j+1,2k+m}-\psi_{j+1,2k+m+1})$. Given a family of intervals $\D$ together with a scale window $j\in[j_0,J]$, define
\[
S_{\D,j}f=\sum_{I\in\D,\ \abs I=2^{-j}}\ip f{\psi_I}\,\Theta_I,\qquad S^{[j_0,J]}_\D=\sum_{j=j_0}^JS_{\D,j}.
\]
For $\psi=h$ and $m=0$, Petermichl's shift is given by $\Theta_I=\frac1{\sqrt2}(h_{I_\ell}-h_{I_r})$.
\end{definition}

\begin{assumption}[Spectral decay]\label{ass:H}
There exist $C,\alpha>0$ satisfying $\abs{\hatt\psi(t)}\le C\min(\abs t,\abs t^{-\alpha})$.
\end{assumption}
\noindent Compact support together with $\int\psi=0$ yields the bound $\abs{\hatt\psi(t)}\le C\abs t$. For Haar ($\alpha=1$), as well as for the Daubechies, symlet and coiflet families employed below, decay at infinity holds [Dau92, Ch.~7].

We denote the number of vanishing moments of $\psi$ by $N\ge1$. We call $\psi$ \emph{smooth} when $\psi\in C^{1,\delta}(\R)$ holds for some $\delta\in(0,1)$.

\paragraph{Which generators are covered.} The H\"older exponent for db2 equals $\alpha\approx0.550$, which means that db2 is not $C^1$; for db4, the corresponding exponent is $\alpha\approx1.619$ [DL92]. As another spectral factor of the same filter polynomial, sym4 shares $\abs{\hatt\psi}$ with db4 and consequently has the same Sobolev exponent $s_0\approx1.78$ [Oja98]. Its H\"older exponent is therefore at least $s_0-\frac12\approx1.28$. In the sense specified above, both db4 and sym4 are thus smooth. By contrast, only the multiplier results of \cref{sec:mult} cover db2: neither \cref{thm:cov} (the $C^1$ case) nor \cref{sec:uniform} applies to it. The H\"older exponents for sym6 and coif3 have not been checked, and we make no use of them; Haar and sym4 are used for all deformation experiments in \cref{sec:num}.

\paragraph{Constants.} Except where otherwise specified, the constant denoted by $C$ depends exclusively on $\psi$, $L$, the child offset $m$ and, in the case of smooth generators, $\delta$. Its value can vary between lines, but it is always independent of $J$, $\beta$, $\tau$ and $f$. Explicit constants are provided in \cref{sec:finiteJ}, whereas the constants furnished by the $T(1)$ theorem in \cref{sec:uniform} are not explicit.

\begin{lemma}[Moments of $\Theta$]\label{lem:mom}
Vanishing moments of orders $0,\dots,N-1$ for $\psi$ imply vanishing moments of orders $0,\dots,N$ for $\Theta^{(m)}$. In the Haar case, $\int s^2\Theta(s)\,ds=\frac1{16}$.
\end{lemma}
\begin{proof}
$\int s^p\Theta(s)\,ds=2^{-p-1}\int[(t+m)^p-(t+m+1)^p]\psi(t)\,dt$, with the bracket being a polynomial in $t$ of degree $p-1$. In the Haar case with $p=2$, the integral equals $-\frac14\int t\,h(t)\,dt=\frac1{16}$.
\end{proof}

\subsection{Random lattices}
\begin{assumption}[Translation model]\label{ass:trans}
$\beta\sim U[0,1)$ and $\D^\beta=\{2^{-j}([0,1)+k)+\beta:\ j\ge0,\ k\in\Z\}$. Our notation is $S^J_\beta:=S^{[0,J]}_{\D^\beta}$ and $S_{\beta,j}:=S_{\D^\beta,j}$.
\end{assumption}
For $j\ge0$, the relation $1\in2^{-j}\Z$ gives $\D^{\beta+1}=\D^\beta$. Since $\beta\sim U[0,1)$, for each $j\ge0$ the offset $\beta\bmod2^{-j}$ is uniform on $[0,2^{-j})$. The coarsest scale is $j=0$, whose cells have length $1$.

\begin{definition}[All-scale random dyadic system; dilations]\label{def:allscale}
Take $\omega=(\omega_i)_{i\in\Z}$ with i.i.d.\ Bernoulli$(\frac12)$ coordinates, and let $x_j(\omega)=\sum_{i>j}\omega_i2^{-i}$ and $\D^\omega=\{2^{-j}([0,1)+k)+x_j(\omega):\ j,k\in\Z\}$, Hyt\"onen's random dyadic system [H12]. For each $j\in\Z$, its offset at scale $j$ is uniform on $[0,2^{-j})$; restricted to $j\ge0$, it is the translation model with $\beta=\sum_{i\ge1}\omega_i2^{-i}$. Given $r\in[1,2)$, set $D_rf(x)=r^{-1/2}f(x/r)$ and $S_{\omega,r}:=D_rS_{\D^\omega}D_r^{-1}$. In the \emph{dilation model}, a draw of $r$ is taken from $\frac{dr}{r\ln2}$ on $[1,2)$ independently of $\omega$.
\end{definition}

\subsection{Feature map and deformations}
Let $\phi\ge0$, $\int\phi=1$, $\phi'\in L^1$, $\phi_J(x)=2^{-J}\phi(2^{-J}x)$ and $A_Jf=f*\phi_J$. Define
\[
\Psi_1'f=\E_\beta\abs{S^J_\beta f},\qquad\Psi_1=A_J\Psi_1',\qquad\Psi_0=A_J,\qquad\Phi f=(\Psi_0f,\Psi_1f).
\]
A deformation takes the form $L_\tau f(x)=f(x-\tau(x))$, where $\tau\in C^1$ and $\eps:=\norm{\tau'}_\infty\le\frac12$. Consequently, $\rho:=\mathrm{id}-\tau$ is a $C^1$ diffeomorphism satisfying $\rho'\in[\frac12,\frac32]$, while $\rho^{-1}$ is $2$-Lipschitz. Throughout all results below we impose this standing assumption, thereby ensuring that $\rho$ remains increasing and invertible. In addition, we employ the mass-preserving transport $\Lambda_\tau f:=\rho'\cdot(f\circ\rho)$, which satisfies $\norm{\Lambda_\tau f}^2=\int\rho'(\rho^{-1}u)\abs{f(u)}^2du\le\frac32\norm f^2$.

\paragraph{Scale-separated channels and a second layer.} Mallat's scattering retains a single channel for each scale. The DSSN channels corresponding to these are
\[
\Psi_{1,j}'f=\E_\beta\abs{S_{\beta,j}f}\qquad(0\le j\le J).
\]
For a collection $\mathcal P\subset\{0,\dots,J\}^2$ of scale pairs (such as $j_2<j_1$), an independent lattice offset $\beta'$ is used in the second layer:
\[
\Psi_{2,(j_1,j_2)}'f=\E_{\beta'}\abs{S_{\beta',j_2}\Psi_{1,j_1}'f},\qquad(j_1,j_2)\in\mathcal P .
\]
The feature map of depth two is $\Phi^{(2)}f=\big(A_Jf,\ (A_J\Psi_{1,j}'f)_{j},\ (A_J\Psi_{2,p}'f)_{p\in\mathcal P}\big)$. To measure families of functions we use $\ell^2(L^2)$, that is, \ $\norm{(u_q)_q}^2=\sum_q\norm{u_q}^2$.

\paragraph{Dilation averaging.} For $r\in[1,2)$ define $S^J_{\beta,r}:=D_rS^J_\beta D_r^{-1}$ as the shift on the dilated lattice $r\D^\beta$. Given a probability measure $\Pi$ on $[1,2)$, define
\[
\Psi_{1,\Pi}'f=\int_{[1,2)}\E_\beta\abs{S^J_{\beta,r}f}\,d\Pi(r).
\]
The dilation model of \cref{def:allscale} is matched by the choice $d\Pi=\frac{dr}{r\ln2}$, whereas $\Pi=\frac1Q\sum_q\delta_{r_q}$ denotes the $Q$-point quadrature employed in practice (\cref{sec:mult}). The summed, translation-averaged channel $\Psi_1'$ is analysed in \Cref{sec:finiteJ} and \cref{sec:uniform}; the results extend to the variants above by \cref{ssec:ext}.

\subsection{Structural facts}
\begin{proposition}[Partial isometry and non-expansiveness]\label{prop:struct}
For each $\psi$, $m$ and $\beta$:
\begin{enumerate}[label=(\alph*),nosep]
\item $\{\Theta_I\}_{I\in\D^\beta,\,0\le j\le J}$ forms an orthonormal system, and $\norm{S^J_\beta f}=\norm{P^\beta_{[0,J]}f}\le\norm f$, with $P^\beta_{[0,J]}$ denoting the orthogonal projection onto $\overline{\operatorname{span}}\{\psi_I:I\in\D^\beta,\,0\le j\le J\}$;
\item $\norm{\Psi_if-\Psi_ig}\le\norm{f-g}$ for $i=0,1$; moreover, $\norm{\Phi f-\Phi g}\le\sqrt2\norm{f-g}$.
\end{enumerate}
\end{proposition}
\begin{proof}
(a) Injectivity of the map $(j,k)\mapsto\{(j+1,2k+m),(j+1,2k+m+1)\}$ and pairwise disjointness of its images imply that the $\Theta_I$ are normalized differences of disjoint pairs belonging to an orthonormal system. (b) $\norm{A_J}\le1$ and $\abs{\abs a-\abs b}\le\abs{a-b}$; next use (a) together with Minkowski's integral inequality.
\end{proof}

\begin{lemma}[Joint measurability]\label{lem:meas}
For $f\in L^2$ and finite $J$, joint measurability holds for $(\beta,x)\mapsto S^J_\beta f(x)$.
\end{lemma}
\begin{proof}
$\beta\mapsto\ip f{\psi_{\beta,j,k}}$ is continuous and $(\beta,x)\mapsto\Theta_{\beta,j,k}(x)$ is Borel; moreover, only finitely many $(j,k)$ contribute on $[0,1)\times[-R,R]$.
\end{proof}

%==================================================================
\section{The averaged multiplier and phase alignment}\label{sec:mult}

This section provides answers to (Q1) and (Q4). Although \Cref{prop:mult,thm:haarkappa} are known in substance [P00, H08], we include brief Fourier-side proofs, as the symbol $m_0$ appears throughout the paper. New here are the quantitative statements (\cref{cor:transfail}, the dilation-quadrature ripple) together with the phase alignment of smooth generators (\cref{ssec:smooth}).

Set
\[
m_0(\xi):=\overline{\hatt\psi(\xi)}\,\hatt\Theta(\xi),\qquad\hatt\Theta(\xi)=\tfrac12\hatt\psi(\xi/2)\,e^{-\pi im\xi}(1-e^{-\pi i\xi}).
\]
By \cref{ass:H} we have $\abs{m_0(t)}\le C'\min(t^2,\abs t^{-2\alpha})$; consequently, $\sum_{j\in\Z}\abs{m_0(2^{-j}\xi)}$ is bounded uniformly in $\xi\ne0$. The fact that $\psi$ is real yields $m_0(-\xi)=\overline{m_0(\xi)}$.

\begin{proposition}[Multiplier of the averaged shift]\label{prop:mult}
\begin{enumerate}[label=(\alph*),nosep]
\item Under \cref{ass:trans}, the Fourier multiplier whose symbol is $m_J(\xi)=\sum_{j=0}^Jm_0(2^{-j}\xi)$ is $\E_\beta S^J_\beta$.
\item In the model of \cref{def:allscale}, the symbol of $\E_\omega S^{[-K,J]}_{\D^\omega}$ is $\sum_{j=-K}^Jm_0(2^{-j}\xi)$. Strong convergence of these operators to the multiplier $m_T(\xi)=\sum_{j\in\Z}m_0(2^{-j}\xi)$ holds as $K,J\to\infty$, and $m_T(2\xi)=m_T(\xi)$.
\item In the dilation model we have
\[
\E_{\omega,r}S_{\omega,r}=\operatorname{Re}\kappa_\psi\cdot\Id-\operatorname{Im}\kappa_\psi\cdot H,\qquad\kappa_\psi:=\frac1{\ln2}\int_0^\infty m_0(t)\,\frac{dt}t .
\]
\end{enumerate}
\end{proposition}
\begin{proof}
(a) Fix $j$, and take $s=2^{-j}$ and $\psi_{s,b}(x)=s^{-1/2}\psi((x-b)/s)$. Uniformity of $\beta\bmod s$ on $[0,s)$ implies $\E_\beta S_{\beta,j}f=\int_\R\ip f{\psi_{s,b}}\Theta_{s,b}\,\frac{db}s$. For Schwartz $f$, substitution of $\hatt{\psi_{s,b}}(\eta)=s^{1/2}\hatt\psi(s\eta)e^{-2\pi ib\eta}$ followed by integration in $b$ yields $m_0(s\xi)\hatt f(\xi)$. (b) At each scale the computation is identical, as it involves only the marginal law of $x_j(\omega)$. Bounded and pointwise convergence of the partial sums therefore gives strong convergence of the operators. (c) If the symbol of $T$ is $m(\xi)$, the symbol of $D_rTD_r^{-1}$ is $m(r\xi)$; consequently, averaging gives an operator with symbol $\frac1{\ln2}\int_1^2m_T(r\xi)\frac{dr}r$. Boundedness of $\sum_j\abs{m_0(2^{-j}\eta)}$ together with $\int_0^\infty\abs{m_0(t)}\frac{dt}t<\infty$ permits the use of Fubini. For $\xi>0$, making the substitution $t=2^{-j}r\xi$ yields
\[
\frac1{\ln2}\int_1^2\sum_{j\in\Z}m_0(2^{-j}r\xi)\frac{dr}r=\frac1{\ln2}\sum_j\int_{2^{-j}\xi}^{2^{-j+1}\xi}m_0(t)\frac{dt}t=\kappa_\psi .
\]
For $\xi<0$ set $\xi=-\abs\xi$. The reality of $\psi$ and $\Theta$ implies $m_0(-t)=\overline{m_0(t)}$; applying the same substitution with $t=2^{-j}r\abs\xi$ then yields $\overline{\kappa_\psi}$. Thus, for $\xi\ne0$ the symbol is $\operatorname{Re}\kappa_\psi+i\operatorname{Im}\kappa_\psi\sgn\xi$. In our convention the symbol of $H$ is $-i\sgn\xi$, so this coincides with the symbol of $\operatorname{Re}\kappa_\psi\cdot\Id-\operatorname{Im}\kappa_\psi\cdot H$.
\end{proof}

\begin{proposition}[Haar constant]\label{thm:haarkappa}
For $\psi=h$ together with $m=0$,
\[
m_0(\xi)=\frac{8i\sin^3(\pi\xi/2)\sin^2(\pi\xi/4)}{\pi^2\xi^2},\qquad\kappa_h=\frac{i\pi}{8\ln2},
\]
and hence $H=c_0\,\E_{\omega,r}S_{\omega,r}$, where $c_0=-8\ln2/\pi\approx-1.76508$.
\end{proposition}
\begin{proof}
Using $\hatt h(\xi)=e^{-\pi i\xi}\,2i\sin^2(\pi\xi/2)/(\pi\xi)$ yields $\hatt\Theta(\xi)=-e^{-\pi i\xi}\,4\sin^2(\pi\xi/4)\sin(\pi\xi/2)/(\pi\xi)$; taking the product then yields $m_0$. For $y=\pi t/2$,
\[
\kappa_h=\frac i{\ln2}\int_0^\infty\frac{\sin^3y\,(1-\cos y)}{y^3}\,dy .
\]
The relation $\int_0^\infty(\sin bx-bx)x^{-3}dx=-\pi b^2/4$ ($b>0$) yields
\begin{equation}\label{eq:sinint}
\int_0^\infty\sum_ka_k\sin(b_ky)\,\frac{dy}{y^3}=-\frac\pi4\sum_ka_kb_k^2\qquad\text{whenever }\sum_ka_kb_k=0 .
\end{equation}
Because $\sin^3y=\frac34\sin y-\frac14\sin3y$ and $\sin^3y\cos y=\frac14\sin2y-\frac18\sin4y$ hold, the value of the integral is $\frac{3\pi}8-\frac\pi4=\frac\pi8$.
\end{proof}
\noindent Equivalently, this is Hyt\"onen's identity [H08]. In his notation, $H=-\frac8\pi\int_1^2\!\int\mathrm{III}_{\beta,r}\,d\mu(\beta)\,\frac{dr}r$ uses the unnormalized measure $\frac{dr}r$, which assigns mass $\ln2$ to $[1,2)$. Under our probability normalization, the same statement becomes $H=-\frac{8\ln2}\pi\E_{\omega,r}S_{\omega,r}=c_0\E_{\omega,r}S_{\omega,r}$. The integral evaluated by quadrature is $0.392699082$, in agreement with $\pi/8$. To check the sign in the Haar case, near $0^+$ we have $\operatorname{Re}m_0=0$ and $\operatorname{Im}m_0>0$; thus $\E_\beta S^J_\beta\sin(2\pi\xi\cdot)=\operatorname{Im}m_J(\xi)\cos(2\pi\xi\cdot)$ for $\xi>0$, whereas $H\sin(2\pi\xi\cdot)=-\cos(2\pi\xi\cdot)$. A discrete implementation recovers both this sign and the values of $\operatorname{Im}m_J$ to three digits.

\begin{corollary}[Translation averaging alone does not give $H$]\label{cor:transfail}
For Haar, nonconstancy of $m_T$ on $(0,\infty)$ follows from the fact that, within every octave, $\operatorname{Im}m_T$ varies over $[0.51029,\,0.62083]$, with $\operatorname{Im}m_T(1)=0.60506$ and $\operatorname{Im}m_T(3/2)=0.51167$, whereas $\operatorname{Re}m_T=0$. It follows that:
\begin{enumerate}[label=(\alph*),nosep]
\item \emph{All-scale model.} No real constant $c$ makes the symbol $c\,m_T+i\sgn$ vanish identically on any octave. Therefore, as $K,J\to\infty$, convergence of $c\,\E_\omega S^{[-K,J]}_{\D^\omega}$ to $H$ fails in the strong operator topology; the limiting error for each $f$ is $\norm{(c\,m_T+i\sgn)\hatt f}$.
\item \emph{Translation model (scales $j\ge0$).} The strong limit of $\E_\beta S^J_\beta$ as $J\to\infty$ is the multiplier $m_+(\xi):=\sum_{j\ge0}m_0(2^{-j}\xi)$ rather than $m_T$. In the case of $\abs\xi\ge\xi_0>0$,
\[
\abs{m_+(\xi)-m_T(\xi)}\le\sum_{k\ge1}\abs{m_0(2^k\xi)}\le\frac{8}{3\pi^2\xi_0^2},
\]
so the conclusion of (a), up to this error, holds for $\hatt f$ with support in $\abs\xi\ge\xi_0$.
\end{enumerate}
\end{corollary}
\begin{proof}
The certification of these values uses interval arithmetic (\cref{app:repro}). With $\varphi(t)=\operatorname{Im}m_0(t)$ taken from \cref{thm:haarkappa}, write $f(s)=\operatorname{Im}m_T(2^s)=\sum_j\varphi(2^{s-j})$. For small $t$, the terms with $2^{s-j}\notin[2^{-21},2^{23}]$ admit the bounds $\abs{\varphi(t)}\le\frac{\pi^3}{16}t^3$ and $\abs{t\varphi'(t)}\le\frac{7\pi^3}{16}t^3$; for large $t$, the corresponding bounds are $\abs{\varphi(t)}\le\frac8{\pi^2t^2}$ and $\abs{t\varphi'(t)}\le\frac{8(2\pi+2/t)}{\pi^2t}$. Using the mean-value form $f(s_{\rm mid})+f'(I)\cdot[-\frac h2,\frac h2]$ to enclose $f$ on every subinterval of $s\in[0,1]$, a branch-and-bound search yields $\min f\in[0.5102896,0.5102897]$, $\max f\in[0.6208294,0.6208295]$, $f(0)=0.6050609$ and $f(\log_2\frac32)=0.5116660$ to all displayed digits. Because $\operatorname{Re}m_0=0$, we have $c\,m_T+i=i(c\operatorname{Im}m_T+1)$ on $(0,\infty)$; moreover, log-periodicity together with the computed range shows that $\operatorname{Im}m_T$ is continuous and not constant on any octave. Plancherel and \cref{prop:mult}(b) therefore imply (a). To obtain (b), note that the partial sums $m_J$ tend boundedly and pointwise to $m_+$, while $\sum_{k\ge1}8/(\pi^24^k\xi^2)=8/(3\pi^2\xi^2)$.
\end{proof}

\begin{remark}[The parent--child structure carries the Hilbert component]\label{rem:sib}
\Cref{prop:mult} is applicable to every lattice operator constructed from a fixed pair of generators. Take the dyadic Hilbert transform of Domelevo--Petermichl [DP22]: its action on the two children of each interval is given by $h_{I_r}\mapsto h_{I_\ell}$, $h_{I_\ell}\mapsto-h_{I_r}$. Measuring frequency in units of the child scale as $t$ gives the translation-averaged scale symbol $i\abs{\hatt h(t)}^2\sin2\pi t$ and hence the dilation-averaged constant $\kappa=\frac i{\ln2}\int_0^\infty\abs{\hatt h(t)}^2\sin(2\pi t)\,\frac{dt}t$. For $y=\pi t$, the integrand is a multiple of $\sin^4(y/2)\sin2y/y^3$, and
\[
\sin^4(y/2)\sin2y=\tfrac18\big(-2\sin y+3\sin2y-2\sin3y+\tfrac12\sin4y\big)
\]
obeys $\sum_ka_kb_k=\sum_ka_kb_k^2=0$. Consequently, \eqref{eq:sinint} yields $\kappa=0$, consistent with [DP22]. If only translation averaging is performed, this shift has a log-periodic all-scale symbol with zero mean; within every octave, its imaginary part ranges over $[-0.247,\,0.250]$. Whereas a sibling rotation remains within one scale, the Hilbert component of Petermichl's shift arises through passage from a parent to its children. This passage is also acted on by the phase alignment of \cref{ssec:smooth}.
\end{remark}

\begin{figure}[tbp]\centering
\includegraphics[width=0.80\linewidth]{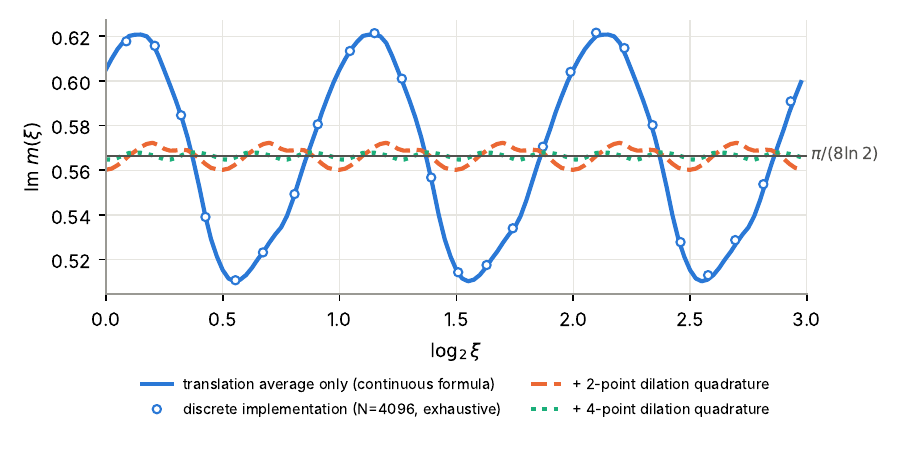}
\caption{Averaged Haar multiplier $\operatorname{Im}m(\xi)$ against $\log_2\xi$. Translation averaging alone ($m_T$, solid; a discrete implementation with $2^{12}$ samples shown as circles) oscillates in $[0.510,0.621]$ with period one octave. $Q$-point dilation quadrature flattens it (\cref{tab:ripple}).}\label{fig:ripple}
\end{figure}

\paragraph{Dilation quadrature.}
In practice we substitute $Q$ log-uniform nodes $r_q=2^{(q+1/2)/Q}$ for $r$. This produces a $\frac1Q$-log-periodic symbol $M_Q(\xi)=\frac1Q\sum_qm_T(r_q\xi)$ whose \emph{ripple} we define to be $\max_\xi\abs{M_Q(\xi)-\kappa}/\abs{\operatorname{Im}\kappa}$: the largest relative deviation from $\kappa$, the dilation-averaged constant and also the mean of $M_Q$ over a log-octave. For Haar, the ripple takes the values $9.930\%$, $1.123\%$, $0.327\%$ and $0.0595\%$ for $Q=1,2,4,5$. Applying the interval method from the proof of \cref{cor:transfail} to $M_Q$ certifies these values to the displayed digits. For $Q=1$, the deviation lies between $-9.93\%$ and $+9.58\%$, giving the introduction's ``about $\pm10\%$''. Refer to \cref{fig:ripple} and \cref{tab:ripple}.

\paragraph{Monte Carlo ensembles.}
For i.i.d.\ $\beta_1,\dots,\beta_M$ and $\hat H_M=\frac{c_0}M\sum_{i=1}^MS^J_{\beta_i}$, combining \cref{prop:struct}(a) with the bias--variance identity yields
\[
\E\norm{\hat H_Mf-Hf}^2=\frac VM+\norm{c_0\E S^Jf-Hf}^2,\qquad V=c_0^2\E\norm{P^\beta_{[0,J]}f}^2-\norm{c_0\E S^Jf}^2\le c_0^2\norm f^2 .
\]
Under averaging over translations alone, the irreducible ripple of \cref{cor:transfail} remains in the bias.

\subsection{Smooth generators and phase alignment}\label{ssec:smooth}
\begin{figure}[tbp]\centering
\includegraphics[width=0.38\linewidth]{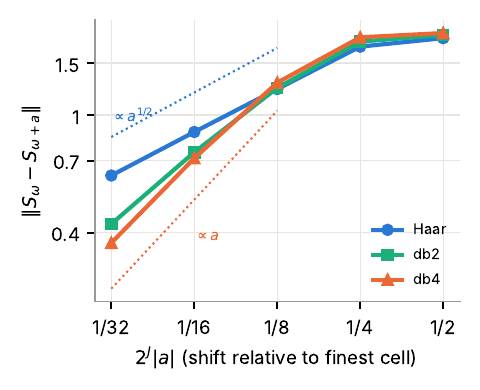}\hfill
\includegraphics[width=0.58\linewidth]{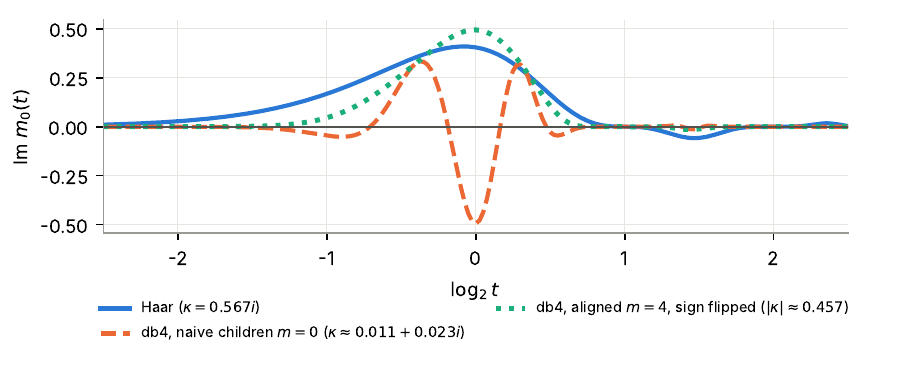}
\caption{Left: single-lattice translation sensitivity $\norm{S^J_\beta-S^J_{\beta+a}}$ (power iteration, $2^{14}$ samples, finest cell 32 samples). Haar behaves like $a^{1/2}$, db4 close to $a$. Right: integrand $\operatorname{Im}m_0(t)$ of the Hilbert constant, whose area is $\ln2\cdot\operatorname{Im}\kappa_\psi$. Haar is single-signed. The naive db4 child pair ($m=0$) winds in phase and cancels, and the aligned pair ($m=4$) restores a single sign.}\label{fig:sample}
\end{figure}

By \cref{prop:mult}(c), the Hilbert coefficient of the averaged shift is $-\operatorname{Im}\kappa_\psi$. This coefficient is obtained by integrating, over relative scales, the cross-spectrum between the parent $\psi$ and the child difference $\Theta$. For Haar, a step-function quadrature pair on the \emph{same} interval is formed by $h=\sgn\sin2\pi x$ and $\Theta=\sgn\cos2\pi x$; consequently, $\operatorname{Im}m_0$ has one sign (\cref{fig:sample}, right). For a Daubechies wavelet supported on $[0,L]$, the parent's support has the standard children $(2k,2k+1)$ at its left end. Because a spatial displacement introduces a winding phase as a multiplier of the cross-spectrum, cancellation occurs in the integral: for db4 the naive ($m=0$) constant $\abs{c_0}=1/\abs{\operatorname{Im}\kappa}$ is $43$, and for coif3 it is $1851$. The two spectra have magnitudes that overlap well; thus phase, rather than a lack of smoothness, causes the cancellation.

\begin{definition}[Aligned generator]
The \emph{aligned} child offset $m^\ast$ is the offset that maximizes $\abs{\operatorname{Im}\kappa_{\psi,m}}$ over $m\in\{0,\dots,L-1\}$.
\end{definition}
\noindent The symbol depends on the offset solely via the factor $e^{-\pi im\xi}$ in $\hatt\Theta$; consequently, $m_0^{(m)}(t)=e^{-\pi imt}m_0^{(0)}(t)$. Thus, selecting $m$ amounts to selecting an integer linear phase that most effectively compensates the phase of $\overline{\hatt\psi(t)}\hatt\psi(t/2)(1-e^{-\pi it})$ over the range of $t$ on which $\abs{m_0}$ is large, rather than setting a free numerical parameter. Since the child map remains injective and $\supp\Theta^{(m)}\subset[0,L]$, all structural statements in this paper apply to every $m$; the regularity hypotheses on $\psi$ are retained in the stability theorems. Only the size of the Hilbert component is affected by alignment.

\begin{table}[htbp]\centering\footnotesize\setlength{\tabcolsep}{4pt}
\begin{tabular}{lcc c c c c c c c}\toprule
gen. & $L$ & $m^\ast$ & $\kappa_{\psi,m^\ast}$ & $\abs{c_0}$ & var. & $\frac{\abs{\operatorname{Re}\kappa}}{\abs{\operatorname{Im}\kappa}}$ & $\abs{c_0}$, $m{=}0$ & $Q{=}1$ & $Q{=}4$\\\midrule
Haar & 1 & 0 & $0.56655i$ & 1.765 & $1\times$ & $0$ & --- & 9.93\% & 0.327\%\\
db2 & 3 & 1 & $0.0424+0.5272i$ & 1.897 & $1.16\times$ & 8.04\% & 3.9 & 16.5\% & 0.145\%\\
db4 & 7 & 4 & $0.0104-0.4566i$ & 2.190 & $1.54\times$ & 2.27\% & 43.2 & 24.4\% & 0.131\%\\
sym4 & 7 & 3 & $-0.0113+0.4610i$ & 2.169 & $1.51\times$ & 2.45\% & 85.7 & 22.4\% & 0.075\%\\
sym6 & 11 & 5 & $0.0013+0.4017i$ & 2.489 & $1.99\times$ & 0.32\% & 226 & 32.6\% & 0.049\%\\
coif3 & 17 & 8 & $0.0064+0.3910i$ & 2.558 & $2.10\times$ & 1.64\% & 1851 & 35.0\% & 0.052\%\\\bottomrule
\end{tabular}
\caption{Aligned generators. Values are computed in double precision from the spectrum $\hatt\psi(\xi)=\mathsf p(\xi/2)\prod_{k\ge1}\mathsf q(\xi/2^{k+1})$, where $\mathsf q$ and $\mathsf p$ are the low- and high-pass filter symbols of the PyWavelets filters, truncated after $40$ factors; no sampled wavelet enters. An independent recomputation (midpoint rule on $(0,600]$ with step $0.004$) reproduces every entry to the digits shown. The ripple columns for smooth generators are double-precision values, not certified bounds. ``var.'' is the Monte Carlo variance relative to Haar, $(c_0/c_{0,h})^2$. ``$\abs{c_0}$, $m{=}0$'' is the unaligned generator. The last two columns are ripples under $Q$-point dilation quadrature. db2 is not $C^1$ and is covered only by \cref{sec:mult} (see \cref{sec:setting}).}\label{tab:aligned}\label{tab:ripple}
\end{table}

\noindent Findings from \cref{tab:aligned}:
\begin{itemize}[nosep,leftmargin=1.2em]
\item $m^\ast=(L-1)/2$ for the near-linear-phase families db2, sym4, sym6 and coif3. For the extremal-phase db4, $m^\ast=4$ lies one child cell past the support midpoint.
\item $\operatorname{Re}\kappa\ne0$ for each aligned smooth generator. Haar alone among compactly supported real orthonormal wavelets has a symmetry [Dau92], so exact reflection antisymmetry is impossible. Subtracting the identity component, $\operatorname{Re}\kappa\cdot\Id$, or averaging over reflections removes it. Of the tested families, sym6 yields the smallest identity component ($0.32\%$), with twice the Haar variance.
\item Under translation averaging alone, the ripple of the aligned generators is \emph{more} than that of Haar (16--35\% versus 9.9\%); dilation quadrature is therefore necessary. Their smoother log-periodic symbols make them ripple \emph{less} with $Q=4$ (0.05--0.15\% versus 0.33\%).
\end{itemize}

\paragraph{Limits of the Haar generator.}
The averaged network requires no smoothness for stability (\cref{sec:finiteJ}). Nevertheless, Haar presents additional limitations that become relevant when DSSN serves as a surrogate. (i) Approximation order saturates, with $\norm{P_{>J}f}^2=\frac1{48}2^{-2J}\norm{f'}^2(1+o(1))$ for $f\in C^2_c$. (ii) Lattice dependence is merely H\"older-$\frac12$: taking $\norm{h_I-h_{I+a}}^2=6\abs a/\abs I$ for $\abs a\le\abs I/2$ yields $\norm{S^J_\beta-S^J_{\beta+a}}\lesssim(2^J\abs a)^{1/2}$, whereas $\psi\in C^1$ yields $\lesssim2^J\abs a$. The measured exponents are $0.49$ for Haar, $0.76$ for db2 and $0.90$ for db4 (\cref{fig:sample}, left). (iii) Because Haar functions are not in $H^{1/2}$, Sobolev estimates for the Haar network extend only to $s<\frac12$. All three limitations disappear with smooth aligned generators. For these generators, \cref{sec:uniform} establishes scale-uniform bounds for general $C^2$ warps; for Haar, such bounds are still conjectural (\cref{conj:haar}).

%==================================================================
\section{Deformation covariance at finite depth}\label{sec:finiteJ}

\begin{remark}[A $\norm\tau_\infty$ term is necessary for $\Phi$]\label{rem:counter}
No estimate of the form $\norm{\Phi(L_\tau f)-\Phi f}\le C\norm{\tau'}_\infty\norm f$ is possible. Indeed, when $\tau\equiv a$, one has $\tau'=0$, whereas $\norm{\Psi_0T_af-\Psi_0f}=\norm{T_ag-g}\to\sqrt2\norm g$ as $\abs a\to\infty$, with $g=A_Jf$. A $2^{-J}\norm\tau_\infty$ term is required in the correct form; see [M12]. In contrast, no term of this kind is needed for the \emph{covariance} error of $\Psi_1'$.
\end{remark}

\subsection{Local translation identity}
\begin{lemma}[Local translation identity]\label{lem:local}
Given any $f\in L^2$, $a\in\R$ and $x$: $S^J_\beta T_af(x)=S^J_{\beta-a}f(x-a)$, which implies $\Psi_1'(T_af)(x)=\Psi_1'f(x-a)$. Thus, whenever $\tau$ is measurable and for any $x$,
\[
(L_\tau\Psi_1'f)(x)=\Psi_1'\big(T_{\tau(x)}f\big)(x).
\]
\end{lemma}
\begin{proof}
$\psi_{\beta,j,k}(y+a)=\psi_{\beta-a,j,k}(y)$ together with $\Theta_{\beta,j,k}(x)=\Theta_{\beta-a,j,k}(x-a)$. The relations $S^J_{\beta+1}=S^J_\beta$ and the fact that $\beta\mapsto\beta-a\bmod1$ preserves the uniform law give $\E_\beta\abs{S^J_{\beta-a}f(x-a)}=\Psi_1'f(x-a)$. To establish the last claim, choose $x$ and set $a=\tau(x)$.
\end{proof}

Combining \cref{lem:local} with $\abs{\abs u-\abs v}\le\abs{u-v}$ gives
\begin{equation}\label{eq:reduce}
\abs{\Psi_1'(L_\tau f)(x)-L_\tau\Psi_1'f(x)}\le\E_\beta\abs{Q_\beta f(x)},\qquad Q_\beta f(x):=S^J_\beta\big(L_\tau f-T_{\tau(x)}f\big)(x),
\end{equation}
and an application of Minkowski's integral inequality (justified by \cref{lem:meas}), followed by Jensen, yields
\begin{equation}\label{eq:reduce2}
\norm{\Psi_1'(L_\tau f)-L_\tau\Psi_1'f}\le\big(\E_\beta\norm{Q_\beta f}^2\big)^{1/2}.
\end{equation}
The translation $\tau(x)$ is absorbed exactly by the lattice average, leaving only the local distortion of $\tau$ within a single wavelet support.

\subsection{Exact reductions}
\begin{proposition}[The Jacobian part costs no $J$]\label{prop:jac}
Given arbitrary $\beta$ and $x$,
\[
Q_\beta f(x)=S^J_\beta\big[\tau'L_\tau f\big](x)+G_\beta f(x),\qquad G_\beta f(x):=S^J_\beta[\Lambda_\tau f](x)-S^J_\beta[T_{\tau(x)}f](x),
\]
and therefore
\begin{equation}\label{eq:jac}
\norm{\Psi_1'(L_\tau f)-L_\tau\Psi_1'f}\le\sqrt2\,\eps\norm f+\big(\E_\beta\norm{G_\beta f}^2\big)^{1/2}.
\end{equation}
The kernel of $G_\beta$ at scale $j$ takes the form $D_{\beta,j}(x,u)=k_{\beta,j}(x,\rho^{-1}u)-k_{\beta,j}(x,u+\tau(x))$, with $k_{\beta,j}$ denoting the kernel of $S_{\beta,j}$.
\end{proposition}
\begin{proof}
$L_\tau f=\Lambda_\tau f+\tau'L_\tau f$, while $S^J_\beta$ is linear. Apply $\norm{S^J_\beta}\le1$ together with $\norm{L_\tau f}^2=\int\abs f^2(\rho^{-1})'\le2\norm f^2$. Substituting $u=\rho(z)$ gives the kernel.
\end{proof}

\begin{proposition}[Orthogonal splitting]\label{prop:split}
For every lattice interval $I$, choose a reference point $c_I$ and define $\varphi_I:=\psi_I\circ\rho^{-1}-\psi_I(\cdot+\tau(c_I))$,
\[
M_\beta f:=\sum_{I:\,0\le j\le J}\ip f{\varphi_I}\Theta_I,\qquad R_\beta f(x):=\sum_{I:\,0\le j\le J}\big\langle f,\psi_I(\cdot+\tau(c_I))-\psi_I(\cdot+\tau(x))\big\rangle\,\Theta_I(x).
\]
Then one obtains $G_\beta=M_\beta+R_\beta$ together with $\norm{M_\beta f}^2=\sum_I\abs{\ip f{\varphi_I}}^2$; orthogonality across all scales is exact.
\end{proposition}
\begin{proof}
The coefficient corresponding to $\Theta_I(x)$ in $G_\beta f(x)$ equals $\ip f{\psi_I\circ\rho^{-1}}-\ip f{\psi_I(\cdot+\tau(x))}$. Add and subtract $\ip f{\psi_I(\cdot+\tau(c_I))}$, then apply \cref{prop:struct}(a).
\end{proof}

\subsection{Scale-wise bounds}
\begin{proposition}[Haar]\label{prop:haar-MR}
Let $\psi=h$ and $m=0$, and let $c_I$ denote the midpoint of $I$. It follows that
\[
\big(\E_\beta\norm{M_\beta f}^2\big)^{1/2}\le2\sqrt{J+1}\,\eps\norm f,\qquad\big(\E_\beta\norm{R_{\beta,j}f}^2\big)^{1/2}\le\sqrt{48}\,\eps\norm f\quad\text{for each }j,
\]
where the part of $R_\beta$ at scale $j$ is denoted by $R_{\beta,j}$.
\end{proposition}
\begin{proof}
\emph{$M_\beta$.} Let $I=[a,b)$ have midpoint $c$. At $\rho(a),\rho(c),\rho(b)$ the function $h_I\circ\rho^{-1}$ jumps, while the jumps of $h_I(\cdot+\tau(c))$ occur at $a-\tau(c),c-\tau(c),b-\tau(c)$. Order is preserved by $\rho$, and the middle jumps coincide. Consequently, on the two strips between $\rho(g)$ and $g-\tau(c)$, $g\in\{a,b\}$, both of width $\abs{\tau(g)-\tau(c)}\le\eps\abs I/2$, we have $\varphi_I=\pm\abs I^{-1/2}$; elsewhere, $\varphi_I=0$. By Cauchy--Schwarz, $\abs{\ip f{\varphi_I}}^2\le\eps\sum_{g\in\{a,b\}}\int_{\mathrm{strip}_g}\abs f^2$. For a point $u$ in the strip of $g$, we have $\abs{\rho^{-1}u-g}\le\eps2^{-j}$. Each $g\in\mathcal Z_j:=2^{-j}\Z+\beta$ serves as an endpoint for two intervals, hence
\[
\sum_{\abs I=2^{-j}}\abs{\ip f{\varphi_I}}^2\le2\eps\int\abs{f(u)}^2\,\#\{g\in\mathcal Z_j:\abs{\rho^{-1}u-g}\le\eps2^{-j}\}\,du .
\]
The uniformity of $\beta\bmod2^{-j}$, together with the independence of $\rho^{-1}u$ from $\beta$, yields the expected count $2\eps$. Using \cref{prop:split}, sum over $j$.

\emph{$R_{\beta,j}$.} For $x\in I$, we have $r_I(x):=\ip f{h_I(\cdot+\tau(c_I))-h_I(\cdot+\tau(x))}$ as well as $\abs{\tau(c_I)-\tau(x)}\le\eps\abs I/2$. Three strips, each of width $\abs{\tau(c_I)-\tau(x)}$, support the difference of the two translates; they lie at the jumps $g$, namely the two ends and the midpoint, and have heights $\abs I^{-1/2}$, $2\abs I^{-1/2}$ and $\abs I^{-1/2}$. Applying Cauchy--Schwarz yields $\abs{r_I(x)}^2\le3\eps\sum_g\int_{\mathrm{strip}_g(x)}\abs f^2$. We then have $\norm{R_{\beta,j}f}^2=\sum_I\abs I^{-1}\int_I\abs{r_I}^2dx$. With $u$ fixed, the set of $x\in I$ satisfying $u\in\mathrm{strip}_g(x)$ has measure bounded above by $\abs I$. This set is empty unless $\abs{u+\tau(c_I)-g}\le\eps\abs I/2$, a condition that forces $\abs{\rho^{-1}u-g}\le2\eps2^{-j}$. No point of $\mathcal Z_{j+1}=2^{-j-1}\Z+\beta$ is a jump of more than two intervals of scale $j$. Thus $\norm{R_{\beta,j}f}^2\le6\eps\int\abs{f(u)}^2\#\{g\in\mathcal Z_{j+1}:\abs{\rho^{-1}u-g}\le2\eps2^{-j}\}\,du$, giving an expected count of $8\eps$.
\end{proof}

\begin{proposition}[$C^1$ generators]\label{prop:C1}
Let $\psi\in C^1$ and $K_1:=(L+1)\norm\Theta_\infty\norm{\psi'}_\infty$. The scale-$j$ part of $G_\beta$ satisfies $\norm{G_{\beta,j}}\le10L^2K_1\eps$ for each $\beta$ and $j$.
\end{proposition}
\begin{proof}
Write $k_{\beta,j}(x,z)=2^j\varkappa(2^j(x-\beta),2^j(z-\beta))$ with $\varkappa(s,t)=\sum_k\Theta(s-k)\psi(t-k)$, vanishing for $\abs{s-t}>L$ and satisfying $\abs{\partial_t\varkappa}\le K_1$. Let $z_1=\rho^{-1}u$ and $z_2=u+\tau(x)$; thus $z_1-z_2=\tau(z_1)-\tau(x)$. Whenever $D_{\beta,j}(x,u)\ne0$, one has either $\abs{x-z_1}\le L2^{-j}$ or $\abs{x-z_2}\le L2^{-j}$, and in either case
\begin{equation}\label{eq:loc}
\abs{x-z_1}\le2L2^{-j},\qquad\abs{x-z_2}\le2L2^{-j},\qquad\abs{z_1-z_2}\le2\eps L2^{-j}.
\end{equation}
So, with $x$ fixed, the measure of the $u$-support is $\le6L2^{-j}$, whereas fixing $u$ gives measure $\le4L2^{-j}$ for the $x$-support. The mean value theorem yields $\abs{D_{\beta,j}}\le2^{2j}K_1\abs{z_1-z_2}\le2LK_1\eps2^j$; applying Schur's test then yields $2LK_1\eps\sqrt{24}\,L\le10L^2K_1\eps$.
\end{proof}

\begin{theorem}[Finite-$J$ deformation covariance]\label{thm:cov}
Let $\psi$ be either Haar or a compactly supported $C^1$ orthonormal wavelet with any child offset, and take $\eps=\norm{\tau'}_\infty\le\frac12$. We then have
\[
\norm{\Psi_1'(L_\tau f)-L_\tau\Psi_1'f}\le
\begin{cases}\big(\sqrt2+2\sqrt{J+1}+\sqrt{48}\,(J+1)\big)\,\eps\norm f&\text{(Haar)},\\[2pt]
\big(\sqrt2+10L^2K_1(J+1)\big)\,\eps\norm f&(\psi\in C^1).\end{cases}
\]
In particular, for every $J$ the Haar constant does not exceed $10.4(J+1)$, and it equals $(\sqrt{48}+o(1))(J+1)$ as $J\to\infty$. No $\norm\tau_\infty$ term appears in the bound.
\end{theorem}
\begin{proof}
Combine \eqref{eq:jac} with either \cref{prop:split,prop:haar-MR} in the Haar case or \cref{prop:C1} in the smooth case, applying Minkowski over scales.
\end{proof}

\begin{proposition}[A single Haar lattice is only $\eps^{1/2}$-stable]\label{prop:single-haar}
Let $\psi=h$. For every $\beta$ and $j$, we have $\norm{G_{\beta,j}}\le\sqrt{864}\,\eps^{1/2}$. The exponent $\frac12$ is sharp: given $j=0$, $\beta=0$ and $\tau(x)=\eps x$, one can find $f$ satisfying $\norm f=1$ and $\norm{G_{0,0}f}\ge\sqrt{(1-\eps)\eps/8}$.
\end{proposition}
\begin{proof}
Write the scale-$j$ kernel of $G_\beta$ as $D=D_{\beta,j}$, and apply \eqref{eq:loc} with $L=1$. In the Haar case, the sole contribution to $\varkappa(s,t)$ comes from $k=\lfloor s\rfloor$; consequently, unless $\lfloor s\rfloor=\lfloor t\rfloor$, we have $\abs\varkappa\le1$ and $\varkappa(s,t)=0$.

\emph{Upper bound.} For fixed $x$, on every cell of $\mathcal Z_{j+1}=2^{-j-1}\Z+\beta$ the function $z\mapsto k_{\beta,j}(x,z)$ is constant. Thus $D(x,u)\ne0$ requires the segment $[z_1,z_2]$ to contain a point of $\mathcal Z_{j+1}$, that is, \ $u\in B_\beta:=\{u:\operatorname{dist}(\rho^{-1}u,\mathcal Z_{j+1})\le2\eps2^{-j}\}$; also, $\abs D\le2\cdot2^j$. When $a=2^{-j-1}$, no more than $9$ cells meet the window $\abs{x-z_1}\le2\cdot2^{-j}$, and each contributes $z_1$-measure $\le8\eps a$. Consequently $\sup_x\int\abs D\,du\le108\eps$, with Jacobian $\le\frac32$, whereas $\sup_u\int\abs D\,dx\le8$. Cauchy–Schwarz now yields $\norm{G_{\beta,j}f}^2\le\big(\sup_x\int\abs D\,du\big)\int_{B_\beta}\abs{f(u)}^2\big(\int\abs D\,dx\big)du\le864\,\eps\norm f^2$, because $D(x,u)=0$ for $u\notin B_\beta$.

\emph{Sharpness.} Take $j=0$, $\beta=0$ and $\tau(x)=\eps x$, which yield $\rho(x)=(1-\eps)x$; choose $X=[0,\frac14)$ and $U=\rho([\frac12,\frac12+\frac\eps8])$ as well, giving $\abs U=(1-\eps)\eps/8$. When $x\in X$, the sole contribution is from $k=0$, and $\Theta(x)=h(2x)=1$. When $u\in U$, we have $z_1=u/(1-\eps)\in[\frac12,\frac12+\frac\eps8]$ and hence $h(z_1)=-1$. Moreover, combining $z_2=z_1-\eps(z_1-x)$ with $z_1-x\ge\frac14$ yields $\frac{1-\eps}2\le z_2\le\frac12-\frac\eps8$, whence $h(z_2)=1$. Consequently, on $X\times U$ we have $D=-2$. When $f=\abs U^{-1/2}\mathbf 1_U$, we have $\abs{G_{0,0}f}\ge2\abs U^{1/2}$ on $X$, implying $\norm{G_{0,0}f}\ge2\abs U^{1/2}\abs X^{1/2}=\abs U^{1/2}$. Cutting off a version of $\tau$ outside $[-2,2]$ gives $\norm{\tau'}_\infty\le2\eps$ while retaining the same order.
\end{proof}
\noindent The lattice average improves $\eps^{1/2}$ to $\eps$: near a fixed point, the expected number of lattice jumps is proportional to $\eps$ (\cref{prop:haar-MR}). On every lattice, smooth generators are $O(\eps)$-stable (\cref{prop:C1}).

\subsection{Pooling and the feature map}
\begin{lemma}[Pooling]\label{lem:pool}
Provided that $\phi\ge0$, $\int\phi=1$, $\phi'\in L^1$, $\tau\in C^1\cap L^\infty$ and $\eps\le\frac12$ hold, we obtain $\norm{A_JL_\tau-A_J}\le3\eps+2\cdot2^{-J}\norm{\phi'}_1\norm\tau_\infty$.
\end{lemma}
\begin{proof}
The kernel associated with $A_JL_\tau-A_J$ takes the form $P_1+P_2$, where $P_1=\phi_J(x-\rho^{-1}v)[(\rho^{-1})'(v)-1]$ and $P_2=\phi_J(x-\rho^{-1}v)-\phi_J(x-v)$. From $\abs{(\rho^{-1})'-1}\le2\eps$ it follows that $\int\abs{P_1}dx\le2\eps$ and $\int\abs{P_1}dv\le3\eps$; Schur therefore yields $\sqrt6\eps$. When $P_2$, using $z=\rho^{-1}v$ gives $\abs{z-v}\le\norm\tau_\infty$ and consequently $\int\abs{P_2}dx\le\norm\tau_\infty\norm{\phi_J'}_1$. In the $v$-integral, make the substitution $v=\rho(z)$ and express the difference as $\int_0^1\phi_J'(x-z+t\tau(z))\tau(z)\,dt$. Because the derivative of the map $z\mapsto z-t\tau(z)$ is $\ge\frac12$, we obtain $\int\abs{P_2}dv\le3\norm\tau_\infty\norm{\phi_J'}_1$. Applying Schur yields $\sqrt3\norm\tau_\infty\norm{\phi'_J}_1$, together with $\norm{\phi_J'}_1=2^{-J}\norm{\phi'}_1$.
\end{proof}

\begin{corollary}[Mallat-type stability at finite depth]\label{cor:finiteJ}
Assuming the hypotheses of \cref{thm:cov,lem:pool} and taking $C_\psi(J+1)$ to be the constant in \cref{thm:cov},
\[
\norm{\Phi(L_\tau f)-\Phi f}\le\sqrt2\Big[\big(C_\psi(J+1)+3\big)\eps+2\cdot2^{-J}\norm{\phi'}_1\norm\tau_\infty\Big]\norm f .
\]
\end{corollary}
\begin{proof}
$\Psi_1L_\tau f-\Psi_1f=A_J(\Psi_1'L_\tau f-L_\tau\Psi_1'f)+(A_JL_\tau-A_J)\Psi_1'f$, together with $\norm{\Psi_1'f}\le\norm f$. To handle $\Psi_0$, \cref{lem:pool} alone suffices; Pythagoras then combines the two channels.
\end{proof}

%==================================================================
\section{Scale-uniform stability: single lattices and lattice averages}\label{sec:uniform}

The triangle inequality over scales produces the factor $J+1$ appearing in \cref{thm:cov}. Here we investigate when this factor can be removed; by \eqref{eq:jac}, this amounts to a question about $\sup_f\E_\beta\norm{G_\beta f}^2/\norm f^2$.

\subsection{Change of variables and Taylor splitting}
\begin{lemma}[Commutator form]\label{lem:cv}
$G_\beta=\tilde G_\beta\Lambda_\tau$, with $\tilde G_\beta g(x)=\int\tilde K(x,v)g(v)\,dv$ satisfying
\[
\tilde K(x,v)=\sum_{j=0}^J\Big[k_j(x,v)-k_j\big(x,v+\Delta(x,v)\big)\Big],\qquad\Delta(x,v):=\tau(x)-\tau(v),
\]
and $k_j=k_{\beta,j}$. The scale-$j$ term is zero unless $\abs{x-v}\le2L2^{-j}$, because $\abs{\Delta(x,v)}\le\eps\abs{x-v}$.
\end{lemma}
\begin{proof}
$\Lambda_\tau^{-1}h=(h\circ\rho^{-1})(\rho^{-1})'$. Substitute $z=v-\tau(v)+\tau(x)$ into $S^J_\beta T_{\tau(x)}\Lambda_\tau^{-1}h(x)=\int k(x,z)h(\rho^{-1}(z-\tau(x)))(\rho^{-1})'(z-\tau(x))\,dz$.
\end{proof}
\noindent The structure here is that of Calder\'on's first commutator: the kernel compares $k$ at $v$ against $k$ at $v$, displaced by $\tau(x)-\tau(v)$. Set $K_S:=\sum_{j\le J}k_j$ and decompose $\tilde K=F-R$ with
\[
F(x,v)=-\Delta(x,v)\,\partial_vK_S(x,v),\qquad R=\sum_jR_j,\quad R_j=k_j(x,v+\Delta)-k_j(x,v)-\Delta\,\partial_2k_j(x,v).
\]
For smooth $\psi$, take $K_1$ as specified in \cref{prop:C1}, and let $K_\delta$ bound the $C^\delta$ seminorm of $t\mapsto\partial_t\varkappa(s,t)$.

\subsection{A single lattice}
\begin{proposition}[First-order part]\label{prop:F}
Let $\psi\in C^{1,\delta}$. For all $\beta$ and $J$, we then have $\norm F_{2\to2}\le C\eps$ with $C=C(\psi,L,m,\delta)$ (not explicit, as it follows from the $T(1)$ theorem).
\end{proposition}
\begin{proof}
We check the hypotheses of the $T(1)$ theorem as formulated in [Gra14, Sec.~4.3]: suppose that $T$ maps Schwartz functions continuously to tempered distributions, that its kernel obeys the size and smoothness estimates stated below with constant $A$ and exponent $\delta$, that $T$ satisfies the weak boundedness property with constant $B$, and that $T1,T^t1\in\BMO$. The conclusion is $\norm T_{2\to2}\le C_\delta(A+B+\norm{T1}_\BMO+\norm{T^t1}_\BMO)$, where $C_\delta$ depends solely on $\delta$, independently of any a priori bound for $\norm T$. Each hypothesis is verified with constants independent of $J$ and $\beta$.

\emph{(0) A priori boundedness and definitions.} For fixed $J$, boundedness of the kernel $F(x,v)$ follows from $\abs{\Delta\,\partial_vK_S}\le\eps\abs{x-v}\sum_{j\le J}2^{2j}K_1$; moreover, the kernel is zero for $\abs{x-v}>2L$. Schur's test therefore gives boundedness of $F$ on $L^2$, while $F1$ and $F^t1$ are pointwise defined, locally bounded functions. The compact support of the kernel in $x-v$ ensures that these functions coincide with the limits $\lim_RF(\phi(\cdot/R))$ defining $T1$ in [Gra14, Sec.~4.3]. Depending on $J$, this a priori bound serves only to ensure that these objects are well defined.

\emph{(i) Kernel estimates.} When $\abs{x-v}\simeq2^{-\mu}$ holds, the contributing scales are restricted to $j\le\mu+\log_2(2L)$. Consequently, for $\abs h\le\abs{x-v}/2$ we have $\abs{\partial_vK_S}\lesssim K_1\abs{x-v}^{-2}$ and $\abs{\partial_vK_S(x,v+h)-\partial_vK_S(x,v)}\lesssim K_\delta\abs h^\delta\abs{x-v}^{-2-\delta}$. The relation $\Theta,\psi\in C^1$ yields $\abs{\partial_x\partial_vK_S}\lesssim\norm{\Theta'}_\infty\norm{\psi'}_\infty\abs{x-v}^{-3}$ for the mixed derivative. Summing geometrically over the contributing scales produces the stated power in every case, uniformly in $J$ and $\beta$. Combining this with $\abs{\Delta(x,v)}\le\eps\abs{x-v}$ and $\abs{\Delta(x,v)-\Delta(x,v')}\le\eps\abs{v-v'}$ yields
\begin{gather*}
\abs{F(x,v)}\lesssim\frac\eps{\abs{x-v}},\\
\abs{F(x,v)-F(x,v')}+\abs{F(v,x)-F(v',x)}\lesssim\frac{\eps\abs{v-v'}^\delta}{\abs{x-v}^{1+\delta}}\qquad(\abs{v-v'}\le\tfrac12\abs{x-v}).
\end{gather*}
For the second variable, set $F(x,v)-F(x,v')=-[\Delta(x,v)-\Delta(x,v')]\partial_vK_S(x,v)-\Delta(x,v')[\partial_vK_S(x,v)-\partial_vK_S(x,v')]$. Because $\abs{v-v'}\le\abs{x-v}$, the first term, $\lesssim\eps\abs{v-v'}\abs{x-v}^{-2}$, has the stated form; 	the second term is bounded in magnitude by $\lesssim\eps\abs{x-v}\cdot\abs{v-v'}^\delta\abs{x-v}^{-2-\delta}$. In the first variable, combining $\Delta(v,x)-\Delta(v',x)=\tau(v)-\tau(v')$ with the mixed-derivative bound yields $\lesssim\eps\abs{v-v'}\abs{x-v}^{-2}$.

\emph{(ii) Representation.} Integration by parts in $v$ yields, for $f\in C^1_c$,
\[
Ff=-S^J_\beta(\tau'f)+\tau\,S^J_\beta(f')-S^J_\beta(\tau f').
\]
The compact support of $K_S(x,\cdot)$ ensures that no boundary terms occur.

\emph{(iii) Weak boundedness.} Take $\phi,\tilde\phi$ supported in an interval $B$ with centre $c$ and radius $r$, subject to $\norm\phi_\infty,\norm{\tilde\phi}_\infty\le1$ and $\norm{\phi'}_\infty\le r^{-1}$; up to a constant, normalized bumps in the sense of [Gra14] meet this condition. Set $\tau S(f')-S(\tau f')=(\tau-\tau(c))S(f')-S((\tau-\tau(c))f')$, then apply $\abs{\tau-\tau(c)}\le\eps r$ to $B$, $\norm{S^J_\beta}\le1$, $\norm\phi_2,\norm{\tilde\phi}_2\le(2r)^{1/2}$ and $\norm{\phi'}_2\le r^{-1}(2r)^{1/2}$. The contribution from each of the three terms in (ii) is at most $2\eps r$, yielding $\abs{\ip{F\phi}{\tilde\phi}}\le6\eps r$.

\emph{(iv) $F1$.} Assertion (ii) yields $F1=-S^J_\beta(\tau')$. Applying (i) to $K_S$ shows that $S^J_\beta$ is a Calder\'on--Zygmund operator whose norm is $\le1$ and whose constants are independent of $J$ and $\beta$. It therefore maps $L^\infty$ to $\BMO$ uniformly [Gra14, Ch.~4], with $\norm{F1}_\BMO\lesssim\norm{\tau'}_\infty=\eps$.

\emph{(v) $F^t1$.} Because $\int K_S(x,v)\,dx=\sum_I\psi_I(v)\int\Theta_I=0$,
\[
F^t1=-\partial_v\big(S^{J*}_\beta\tau\big)=-\sum_I\ip\tau{\Theta_I}\psi_I'=\sum_I\ip{\tau'}{\tilde\Theta_I}(\psi')_I .
\]
Here $\tilde\Theta(s)=\int_{-\infty}^s\Theta$ has compact support, and $C^1$; by \cref{lem:mom}, $\int\tilde\Theta=-\int s\,\Theta=0$. The $L^2$-normalized lattice copies are the functions $\tilde\Theta_I$ and $(\psi')_I$; moreover, $\ip\tau{\Theta_I}=-2^{-j}\ip{\tau'}{\tilde\Theta_I}$ whereas $\psi_I'=2^j(\psi')_I$. Since $\{\tilde\Theta_I\}$ and $\{(\psi')_I\}$ are dyadic families of mean-zero, compactly supported $C^\delta$ molecules, they are Bessel sequences whose bounds are independent of $J$ and $\beta$. Thus $g\mapsto\sum_I\ip g{\tilde\Theta_I}(\psi')_I$ is uniformly bounded on $L^2$. Summation over the contributing scales once again shows that its kernel $\sum_I\tilde\Theta_I(v)(\psi')_I(x)$ obeys the Calder\'on--Zygmund estimates with exponent $\delta$. It therefore maps $L^\infty\to\BMO$, and $\norm{F^t1}_\BMO\lesssim\eps$.

\emph{(vi) Conclusion.} The conclusion $\norm F_{2\to2}\le C_\delta(A+B+\norm{F1}_\BMO+\norm{F^t1}_\BMO)\le C\eps$ follows from the $T(1)$ theorem, using $A\lesssim\eps$ supplied by (i) and $B\le6\eps$ supplied by (iii). Only $\psi$, $L$, $m$ and $\delta$ determine the constants; the a priori bound in (0) is not used.
\end{proof}

\begin{lemma}[Per-scale remainder]\label{lem:Rj}
Let $\psi\in C^{1,\delta}$. For every $\beta$ and $j$, it follows that $\norm{R_j}_{2\to2}\le C\eps^{1+\delta}$.
\end{lemma}
\begin{proof}
$R_j=\int_0^1[\partial_2k_j(x,v+\theta\Delta)-\partial_2k_j(x,v)]\Delta\,d\theta$; hence, on $\abs{x-v}\le2L2^{-j}$, $\abs{R_j}\le2^{j(2+\delta)}K_\delta\abs\Delta^{1+\delta}\le K_\delta(2L)^{1+\delta}\eps^{1+\delta}2^j$. Schur's test now applies.
\end{proof}

\begin{theorem}[Single lattice]\label{thm:A}
Let $\psi\in C^{1,\delta}$ and $\eps\le\frac12$. For every $\beta$ and $J$, we have $\norm{G_\beta}\le C(\eps+(J+1)\eps^{1+\delta})$, where $C=C(\psi,L,m,\delta)$. The covariance error is therefore bounded above by $C(\eps+(J+1)\eps^{1+\delta})\norm f$; this bound is $J$-uniform provided that $J\lesssim\eps^{-\delta}$.
\end{theorem}
\begin{proof}
Combine \cref{lem:cv}, \cref{prop:F}, \cref{lem:Rj} and \eqref{eq:jac}.
\end{proof}

\begin{lemma}[Kernel of the defect]\label{lem:czk}
Let $\psi\in C^{1,\delta}$. One can choose $A=A(\psi,L,m,\delta)$ so that, for every $J$, every $\beta$ and every $x\ne v$,
\begin{gather*}
\abs{\tilde K(x,v)}\le\frac{A\eps}{\abs{x-v}},\\
\abs{\tilde K(x,v)-\tilde K(x,v')}+\abs{\tilde K(v,x)-\tilde K(v',x)}\le\frac{A\eps\abs{v-v'}^\delta}{\abs{x-v}^{1+\delta}}\qquad\big(\abs{v-v'}\le\tfrac12\abs{x-v}\big).
\end{gather*}
\end{lemma}
\begin{proof}
Write $\tilde K=\sum_j\tilde K_j$ with $\tilde K_j(x,v)=-\int_0^1\partial_2k_j\big(x,v+\theta\Delta(x,v)\big)\Delta(x,v)\,d\theta$, vanishing unless $\abs{x-v}\le2L2^{-j}$. From \emph{Size.} $\abs{\partial_2k_j}\le2^{2j}K_1$ together with $\abs\Delta\le\eps\abs{x-v}\le2L\eps2^{-j}$ we obtain $\abs{\tilde K_j}\le2LK_1\eps2^j$; a sum over scales satisfying $2^j\le2L/\abs{x-v}$ then yields $\abs{\tilde K}\le8L^2K_1\eps/\abs{x-v}$. \emph{Second variable.} By the $C^\delta$ bound for $\partial_t\varkappa$ and $\abs{\Delta(x,v)-\Delta(x,v')}\le\eps\abs{v-v'}$, the integrands evaluated at $v$ and $v'$ differ by at most $2^{j(2+\delta)}K_\delta(2\abs{v-v'})^\delta\cdot\eps\abs{x-v}+2^{2j}K_1\cdot\eps\abs{v-v'}$. The contributing scales all satisfy $2^j\lesssim L/\abs{x-v}$, so their sum yields $\lesssim\eps\abs{v-v'}^\delta\abs{x-v}^{-1-\delta}+\eps\abs{v-v'}\abs{x-v}^{-2}$. This has the stated form because $\abs{v-v'}\le\abs{x-v}$. \emph{First variable.} On the support, $\partial_x\tilde K_j=[\partial_1k_j(x,v)-\partial_1k_j(x,v+\Delta)]-\partial_2k_j(x,v+\Delta)\tau'(x)$ together with $\Theta,\psi\in C^1$ yields $\abs{\partial_x\tilde K_j}\lesssim\eps2^{2j}$. The mean value theorem, followed by summation, provides the bound $\lesssim\eps\abs{v-v'}\abs{x-v}^{-2}$.
\end{proof}

\begin{theorem}[Single-lattice lower bound for affine warps]\label{thm:single-lower}
Let $\psi\in C^{1,\delta}$ for an arbitrary child offset $m$, and take $\tau(x)=c+\eps x$ satisfying $0<\eps\le\frac12$. Set $\lambda=\frac1{1-\eps}$,
\[
c_n(\eps):=\lambda\,\hatt\psi(\lambda n)\,\hatt\Theta(-\eps\lambda n)\quad(n\in\Z),\qquad\sigma_\eps^2:=\sum_{o\in2\Z+1}\Big|\sum_{i\ge0}c_{2^io}(\eps)\Big|^2 .
\]
Then the following hold for each $\beta$:
\begin{enumerate}[label=(\alph*),nosep]
\item $\tilde G_\beta^*1(v)=-\sum_{j=0}^J\Phi_\eps\big(2^j(v-\beta)\big)$, with $\Phi_\eps(t)=\sum_{n\ne0}c_n(\eps)e^{2\pi int}$ having mean zero and being $1$-periodic;
\item $\displaystyle\liminf_{J\to\infty}\frac{\norm{\tilde G_\beta^*1}_{L^2(\beta,\beta+1)}}{\sqrt{J+1}}\ge\sigma_\eps$;
\item there exist $C=C(\psi,L,m,\delta)$ and an absolute constant $c>0$ for which $\norm{G_\beta}\ge c\,\sigma_\eps\sqrt{J+1}-C\eps$ holds for every $J\ge J_0(\eps)$.
\end{enumerate}
In particular, $\sigma_\eps>0$ implies $\sup_J\norm{G_\beta}=\infty$ for every lattice. No assertion of $\sigma_\eps>0$ is made by the theorem; consult \cref{rem:sharp}.
\end{theorem}
\begin{proof}
(a) The compact support of $\tilde K(\cdot,v)$ implies $\tilde G_\beta^*1(v)=\int\tilde K(x,v)\,dx$. 	Express the kernel as $k_j(x,z)=2^j\varkappa(2^j(x-\beta),2^j(z-\beta))$. Integration of the first term yields $\sum_k\psi(2^j(v-\beta)-k)\int\Theta=0$. For the second term, $z=v+\eps(x-v)$; making the substitution $s=2^j(x-\beta)$ yields $\int k_j(x,z)\,dx=\Phi_\eps(2^j(v-\beta))$ with
\[
\Phi_\eps(t)=\int\varkappa\big(s,(1-\eps)t+\eps s\big)\,ds=\sum_{k\in\Z}g_\eps(t-k),\qquad g_\eps(t)=\int\Theta(u)\,\psi\big((1-\eps)t+\eps u\big)\,du ,
\]
following the variable substitution $s=u+k$. Thus the compactly supported function $g_\eps$ has $1$-periodization $\Phi_\eps$, whose Fourier coefficients are given by $\hatt g_\eps(n)$. Using the substitution $w=(1-\eps)t+\eps u$ yields $\hatt g_\eps(\xi)=\lambda\hatt\psi(\lambda\xi)\hatt\Theta(-\eps\lambda\xi)$; consequently, $\hatt g_\eps(n)=c_n$ and $c_0=0$.

(b) Translation allows us to assume $\beta=0$; set $S_J=\sum_{j=0}^J\Phi_\eps(2^j\cdot)$. Express a frequency $N\ne0$ in the form $N=2^Ko$, where $o$ is odd and $K\ge0$. In $S_J$, the coefficient of $e^{2\pi iNt}$ equals $\sum_{j=0}^{\min(J,K)}c_{2^{K-j}o}$. When $K\le J$, this coefficient is the partial sum $P_K(o):=\sum_{i=0}^Kc_{2^io}$. Discarding those frequencies for which $K>J$,
\[
\norm{S_J}_{L^2(0,1)}^2\ge\sum_{o\in2\Z+1}\ \sum_{K=0}^J\abs{P_K(o)}^2 .
\]
Because $\hatt\Theta$ is bounded and $\abs{\hatt\psi(s)}\le C\abs s^{-1-\delta}$, we obtain $\sum_i\abs{c_{2^io}}<\infty$ and therefore $P_K(o)\to B_o:=\sum_{i\ge0}c_{2^io}$. It follows that $\frac1{J+1}\sum_{K\le J}\abs{P_K(o)}^2\to\abs{B_o}^2$ for every $o$; applying Fatou's lemma in $o$ then yields $\liminf_J\norm{S_J}^2/(J+1)\ge\sum_o\abs{B_o}^2=\sigma_\eps^2$. Convergence of the series follows from $\abs{B_o}\le\sum_i\abs{c_{2^io}}\le C\abs o^{-1-\delta}$, which is a consequence of $\abs{c_n}\le C\abs n^{-1-\delta}$.

(c) The mean of each $\Phi_\eps(2^j\cdot)$ on $I=[\beta,\beta+1)$ is zero. Consequently, $\norm{S_J}_{L^2(0,1)}$ is the $L^2$ mean oscillation of $\tilde G_\beta^*1$ on $I$; the John--Nirenberg inequality then yields $\norm{\tilde G_\beta^*1}_\BMO\ge c_0\norm{S_J}_{L^2(0,1)}$. On the other hand, if the kernel of a bounded operator satisfies the estimates of \cref{lem:czk} with constant $A\eps$, both that operator and its adjoint map $L^\infty$ to $\BMO$ with norm at most $C_1\norm{\tilde G_\beta}+C_\delta A\eps$. Here $C_1$ is absolute, whereas $C_\delta$ depends only on $\delta$ [Gra14, Ch.~4]. We therefore obtain $\norm{\tilde G_\beta}\ge c_0C_1^{-1}\norm{S_J}_{L^2(0,1)}-C_1^{-1}C_\delta A\eps$. Lastly, since $\tilde G_\beta=G_\beta\Lambda_\tau^{-1}$ and $\Lambda_\tau$ is a dilation with $\norm{\Lambda_\tau^{-1}}=(1-\eps)^{-1/2}\le\sqrt2$, we have $\norm{G_\beta}\ge\norm{\tilde G_\beta}/\sqrt2$. Combining this with (b) gives the conclusion.
\end{proof}

\begin{remark}[Size of $\sigma_\eps$; Haar]\label{rem:sharp}
$\sigma_\eps$ is determined solely by $\psi$, $m$ and $\eps$ and is independent of $\beta$. Its computation uses double precision and the infinite-product spectrum rather than sampled wavelets; the resulting values are not interval-certified. For aligned sym4 ($m=3$), the values at $\eps=0.05$, $0.1$, $0.2$ are $\sigma_\eps=2.49\times10^{-3}$, $6.49\times10^{-3}$, $3.50\times10^{-2}$, respectively; for db4 ($m=4$), they are $2.62\times10^{-3}$, $6.68\times10^{-3}$, $3.69\times10^{-2}$. Already $B_1\ne0$ in every case. The direct spatial computation of \cref{ssec:num-single} is also reproduced by the exact formula for $\norm{S_J}_{L^2(0,1)}$: at $\eps=0.1$ for sym4, the formula yields $\norm{S_J}/\sqrt{J+1}=(0.164$--$0.169)\,\eps^{1.4}$ for $J=0,\dots,10$, compared with the measured $(0.165\pm0.002)\,\eps^{1.4}$. For each $\eps\in(0,\frac12]$ we expect $\sigma_\eps>0$, with $\sigma_\eps\to0$ faster than $\eps$; a proof, however, would require precise asymptotics of $\hatt\psi$. The equality $\sigma_\eps=0$ would require $B_o(\eps)=0$ for each odd $o$. For Haar, a single lattice does not have a Calder\'on--Zygmund kernel, making \cref{thm:single-lower} inapplicable; instead, the stronger per-scale phenomenon of \cref{prop:single-haar} holds there. Whenever $\sigma_\eps>0$, as indicated for sym4 and db4 by the double-precision values above, a lattice-by-lattice proof of $J$-uniformity is impossible: it must arise from averaging over $\beta$.
\end{remark}

\subsection{Lattice averages are Fourier multipliers}
\begin{lemma}[Averaging produces a multiplier]\label{lem:avg}
Let $T_\beta g(x)=\int K_\beta(x,v)g(v)\,dv$ satisfy $\sup_{\beta,x}\int\abs{K_\beta(x,v)}\,dv<\infty$, with $T_\beta$ bounded on $L^2$ uniformly in $\beta$; assume also $K_{\beta+t}(x+t,v+t)=K_\beta(x,v)$ for every $t$ and $K_{\beta+1}=K_\beta$. Define
\[
\sigma_\beta(\xi,x):=\int K_\beta(x,v)\,e^{2\pi i\xi(v-x)}\,dv ,
\]
which, as a function of $x$, is bounded and $1$-periodic, and satisfies $\sigma_\beta(\xi,x)=\sigma_0(\xi,x-\beta)$. It follows that
\[
\E_\beta\norm{T_\beta g}^2=\int_\R\M(\xi)\abs{\hatt g(\xi)}^2d\xi,\qquad\M(\xi)=\int_0^1\abs{\sigma_0(\xi,u)}^2du\le\sup_u\abs{\sigma_0(\xi,u)}^2 .
\]
\end{lemma}
\begin{proof}
Let $g$ be a Schwartz function. Fubini gives $T_\beta g(x)=\int\hatt g(\xi)e^{2\pi i\xi x}\sigma_0(\xi,x-\beta)\,d\xi$. Choose $\beta$ and express $x=\beta+u+n$ with $u\in[0,1)$ and $n\in\Z$. Periodicity yields $\sigma_0(\xi,x-\beta)=\sigma_0(\xi,u)$; hence
\[
\norm{T_\beta g}^2=\sum_{n\in\Z}\int_0^1\big|F_u(u+n+\beta)\big|^2du,\qquad F_u(t):=\int\hatt g(\xi)\,\sigma_0(\xi,u)\,e^{2\pi i\xi t}\,d\xi .
\]
Integrate with respect to $\beta\in[0,1)$, applying Tonelli. When $n$ ranges over $\Z$ while $\beta$ ranges over $[0,1)$, every point of $\R$ is attained by $t=u+n+\beta$ exactly once. Consequently,
\[
\E_\beta\norm{T_\beta g}^2=\int_0^1\!\!\int_\R\abs{F_u(t)}^2\,dt\,du=\int_0^1\!\!\int_\R\abs{\hatt g(\xi)}^2\abs{\sigma_0(\xi,u)}^2\,d\xi\,du
\]
by Plancherel. Since each side is a bounded quadratic form in $g$, density extends the identity to $L^2$.
\end{proof}

\begin{lemma}[Symbol of a sheared scale-$j$ kernel]\label{lem:sym}
For $\abs b<1$, the symbol of the kernel $K^{(b)}_j(x,v)=k_j(x,v+b(x-v))$ is $\mu^{(b)}(2^{-j}\xi,\,2^j(x-\beta))$, where we take $\lambda=\frac1{1-b}$ and $F_y(s)=\hatt\psi(-s)e^{-2\pi isy}$,
\[
\mu^{(b)}(s,t)=\lambda\sum_k\Theta(t-k)\,F_{t-k}(\lambda s).
\]
\end{lemma}
\begin{proof}
Substitute $w=v+b(x-v)$ to obtain $v-x=\lambda(w-x)$ and $dv=\lambda\,dw$; then apply $\Theta_I(x)\int\psi_I(w)e^{2\pi i\eta(w-x)}dw=\Theta(y)\hatt\psi(-2^{-j}\eta)e^{-2\pi i2^{-j}\eta y}$ with $y=2^j(x-\beta)-k$.
\end{proof}

\begin{lemma}[Spectral bounds]\label{lem:spec}
Let $\Gamma_y(s)=F_y(s)+sF_y'(s)$ and $y\in[0,L]$. Whenever $\psi\in C^{1,\delta}_c$ has $N$ vanishing moments, we have, for $p=0,1$,
\[
\abs{\partial_s^pF_y(s)}\le C\min(\abs s^{N-p},\abs s^{-1-\delta}),\qquad\abs{\partial_s^p\Gamma_y(s)}\le C\min(\abs s^{N-p},\abs s^{-\delta}).
\]
For Haar ($N=1$), these bounds remain valid with $\delta=0$.
\end{lemma}
\begin{proof}
$\hatt\psi$ is entire and has a zero of order $N$ at $0$. When $s$ is large, apply $\abs{\hatt g(s)}\le C\abs s^{-1-\delta}$ for $g\in C^{1,\delta}_c$ to $g=\psi$, $x\psi$ and $x^2\psi$. In the Haar case, bounded variation of these functions yields $\abs s^{-1}$.
\end{proof}

\subsection{Affine warps}
\begin{theorem}[Affine warps, lattice-averaged]\label{thm:B}
Let $\tau(x)=c+ax$, where $\abs a=\eps\le\frac12$. The following holds uniformly in $J$:
\[
\sup_{f\ne0}\frac{\E_\beta\norm{G_\beta f}^2}{\norm f^2}\le\begin{cases}C\eps^2&\psi\in C^{1,\delta}\ (\text{any }m;\ C=C(\psi,L,m,\delta)),\\ C\eps^2\big(1+\log_2\frac1\eps\big)^2&\psi=\text{Haar}\ (C\text{ absolute}).\end{cases}
\]
Hence, independently of $J$, the covariance error for affine warps does not exceed $C\eps\norm f$ (smooth) or $C\eps\log\frac e\eps\norm f$ (Haar).
\end{theorem}
\begin{proof}
For affine $\tau$ we have $\Delta=a(x-v)$; consequently, $\tilde K_\beta$ is covariant in the sense of \cref{lem:avg}. Applying \cref{lem:sym} with $b=0$ and $b=a$ yields the scale-$j$ symbol $\mu(2^{-j}\xi,\cdot)$, with $\mu(s,t)=\sum_k\Theta(t-k)[F_{t-k}(s)-\lambda F_{t-k}(\lambda s)]$. Together, $F(s)-\lambda F(\lambda s)=-\int_1^\lambda\Gamma(\theta s)\,d\theta$ and $\abs{\lambda-1}\le2\eps$ allow \cref{lem:spec} to give
\begin{equation}\label{eq:mubar}
\bar\mu(s):=\sup_t\abs{\mu(s,t)}\le C\begin{cases}\eps\abs s^N&\abs s\le1,\\ \eps\abs s^{-\delta}&1\le\abs s\le\eps^{-1},\\ \abs s^{-1-\delta}\le\eps\abs s^{-\delta}&\abs s\ge\eps^{-1},\end{cases}
\end{equation}
where we estimate $\abs F+\lambda\abs{F(\lambda\cdot)}$ directly in the final range. Thus, when $\delta>0$, $\sum_{j\in\Z}\bar\mu(2^{-j}\xi)\le C\eps$ holds uniformly in $\xi$. For Haar ($\delta=0$), about $\log_2\frac1\eps+1$ dyadic values of $s=2^{-j}\xi$ lie in the middle range $1\le\abs s\le\eps^{-1}$, and each contributes at most $C\eps$. This alone produces the logarithm, whereas the remaining two ranges are geometric. For fixed $J$, no interchange of limits is required: the symbol of $\tilde G_\beta$ is the finite sum of the scale-$j$ symbols. Moreover, by \cref{lem:avg}, $\M(\xi)\le\sup_u\abs{\sum_{j=0}^J\mu(2^{-j}\xi,2^ju)}^2\le(\sum_j\bar\mu(2^{-j}\xi))^2$, and $\Lambda_\tau$ is a dilation with $\norm{\Lambda_\tau}\le\sqrt{3/2}$.
\end{proof}

\noindent Computation from \cref{lem:avg,lem:sym} gives the exact averaged constants $\sup_\xi\M_J(\xi)^{1/2}/\eps$: for sym4 they are $2.58$--$3.22$ at every tested $J$, whereas for Haar they rise with $J$, reaching saturation at $2.37$--$2.52$ by $J\approx8$ (\cref{tab:aff}). With respect to $\eps$, the Haar constant changes very little, suggesting that crude summation introduces the logarithm in \cref{thm:B} as an artefact.

\subsection{$C^2$ warps}\label{ssec:C2}
Let $\curv=\norm{\tau''}_\infty$, choosing units so that the length of the coarsest cell is $1$. Set $\Delta_\ell(x,v)=\tau'(v)(x-v)$; hence $\abs{\Delta-\Delta_\ell}\le\frac{\curv}2\abs{x-v}^2$. Define $R^\ell_j$ by replacing $\Delta$ with $\Delta_\ell$ in $R_j$.

\begin{lemma}[Freezing the slope]\label{lem:freeze}
Let $\psi\in C^{1,\delta}$. We have $\sum_{j\ge0}\norm{R_j-R^\ell_j}_{2\to2}\le C(\eps\curv^\delta+\eps^\delta\curv)$ for each $\beta$.
\end{lemma}
\begin{proof}
Write each in the form $\int_0^1[\partial_2k_j(x,v+\theta D)-\partial_2k_j(x,v)]D\,d\theta$, taking $D=\Delta$ and $D=\Delta_\ell$. On $\abs{x-v}\le2L2^{-j}$, the difference is bounded by $2^{j(2+\delta)}K_\delta(\abs{\Delta-\Delta_\ell}^\delta\abs\Delta+\abs{\Delta_\ell}^\delta\abs{\Delta-\Delta_\ell})\lesssim\eps\curv^\delta2^{j(1-\delta)}+\eps^\delta\curv$. Applying Schur's test introduces a factor of $\lesssim2^{-j}$, summable over $j\ge0$.
\end{proof}

\begin{lemma}[Superposition over slopes]\label{lem:super}
Let $D_{b,j}(x,v)=(x-v)[\partial_2k_j(x,v+b(x-v))-\partial_2k_j(x,v)]$ and $D_b=\sum_jD_{b,j}$. Define $\chi_b(v)=\sgn b$ when $b$ is between $0$ and $\tau'(v)$, and $\chi_b(v)=0$ otherwise. It follows that $R^\ell=\sum_jR^\ell_j$ satisfies
\[
R^\ell g=\int_{-\eps}^\eps D_b(\chi_bg)\,db,\qquad\big(\E_\beta\norm{D_bh}^2\big)^{1/2}\le C\abs b^\delta\norm h,\qquad\big(\E_\beta\norm{R^\ell g}^2\big)^{1/2}\le C\eps^{1+\delta}\norm g,
\]
for $\psi\in C^{1,\delta}$, uniformly in $J$.
\end{lemma}
\begin{proof}
$R^\ell_j=\int_0^{\tau'(v)}\partial_b[k_j(x,v+b(x-v))]\,db-\tau'(v)(x-v)\partial_2k_j(x,v)=\int_0^{\tau'(v)}D_{b,j}\,db$. Lattice covariance holds for $D_b$. Differentiation of \cref{lem:sym} with respect to $b$ gives $\sum_k\Theta(t-k)[\lambda^2\Gamma_{t-k}(\lambda s)-\Gamma_{t-k}(s)]$ as its scale-$j$ symbol, with $s=2^{-j}\xi$. Using \cref{lem:spec},
\begin{align}
\abs{\lambda^2\Gamma(\lambda s)-\Gamma(s)}&\le\abs{\lambda^2-1}\abs{\Gamma(\lambda s)}+\abs{\Gamma(\lambda s)-\Gamma(s)}\nonumber\\
&\lesssim\abs b\min(\abs s^N,\abs s^{-\delta})+\min\big(\abs b\abs s\sup\abs{\Gamma'},\,\abs s^{-\delta}\big),\label{eq:Dbsym}
\end{align}
which is $\lesssim\abs b\abs s^N$ when $\abs s\le1$ and $\lesssim\min(\abs b\abs s^{1-\delta},\abs s^{-\delta})$ when $\abs s\ge1$. Dividing the dyadic sum at $\abs s=\abs b^{-1}$ yields $\lesssim\abs b^\delta$, uniformly with respect to $\xi$ and $J$. For fixed $J$, the scale sum is finite; the second claim follows from \cref{lem:avg}. To obtain the third, apply Minkowski in $b$, together with $\norm{\chi_bg}\le\norm g$ and $\int_{-\eps}^\eps\abs b^\delta db\lesssim\eps^{1+\delta}$.
\end{proof}

\begin{theorem}[$C^2$ warps, lattice-averaged]\label{thm:C}
Let $\psi\in C^{1,\delta}$ for an arbitrary child offset, and $\tau\in C^2$ with $\eps\le\frac12$. With $C=C(\psi,L,m,\delta)$, uniformly over $J$,
\[
\big(\E_\beta\norm{G_\beta f}^2\big)^{1/2}\le C\big(\eps+\eps\curv^\delta+\eps^\delta\curv\big)\norm f,
\]
and consequently $\norm{\Psi_1'(L_\tau f)-L_\tau\Psi_1'f}\le C_\psi(\eps+\eps\curv^\delta+\eps^\delta\curv)\norm f$.
\end{theorem}
\begin{proof}
$\tilde G_\beta=F-(R-R^\ell)-R^\ell$. For each $\beta$, apply \cref{prop:F} and \cref{lem:freeze}; use \cref{lem:super} on average, Minkowski in $L^2(d\beta;L^2)$, and \cref{lem:cv} and \eqref{eq:jac}.
\end{proof}

\begin{theorem}[Scale-uniform stability of the feature map]\label{thm:main}
Take $\psi\in C^{1,\delta}$ with an arbitrary child offset, $\tau\in C^2\cap L^\infty$ satisfying $\eps=\norm{\tau'}_\infty\le\frac12$ and $\curv=\norm{\tau''}_\infty$, and $\phi$ as specified in \cref{lem:pool}. It follows that
\[
\norm{\Phi(L_\tau f)-\Phi f}\le C_\psi\Big[\eps+\eps\curv^\delta+\eps^\delta\curv+2^{-J}\norm{\phi'}_1\norm\tau_\infty\Big]\norm f
\]
where $C_\psi$ does not depend on $J$.
\end{theorem}
\begin{proof}
Argue as in \cref{cor:finiteJ}, replacing \cref{thm:cov} by \cref{thm:C}. If the constant in \cref{thm:C} is $C$, then $\sqrt2\big[C(\eps+\eps\curv^\delta+\eps^\delta\curv)+3\eps+2\cdot2^{-J}\norm{\phi'}_1\norm\tau_\infty\big]\norm f$ bounds the left side from above; hence $C_\psi=\sqrt2\max(C+3,2)$ suffices. Independence from $J$ follows from that of $C$.
\end{proof}

\begin{remark}[Comparison with wavelet scattering]
The structure in $\eps+(\text{curvature terms})+2^{-J}\norm\tau_\infty$ is parallel to Mallat's bound $\norm{\tau'}_\infty\log(\cdot)+\norm{\tau''}_\infty+2^{-J}\norm\tau_\infty$ for wavelet scattering [M12], but DSSN has two specific features. First, since the local translation identity removes translations exactly, no logarithm of $\norm\tau_\infty/\norm{\tau'}_\infty$ occurs in the first-order term. Second, one cannot omit the lattice average: \cref{thm:single-lower} rules out a $J$-uniform bound for affine warps on any single lattice whenever $\sigma_\eps>0$, as the computations for sym4 and db4 in \cref{rem:sharp} show.
\end{remark}

\subsection{Scale-separated channels, a second layer, and dilation averaging}\label{ssec:ext}
Denote the pooling constant of \cref{lem:pool} by $a_J(\tau):=3\eps+2\cdot2^{-J}\norm{\phi'}_1\norm\tau_\infty$.

\begin{proposition}[Non-expansiveness]\label{prop:nonexp}
For every $f,g\in L^2$, we have
$\sum_j\norm{\Psi_{1,j}'f-\Psi_{1,j}'g}^2\le\norm{f-g}^2$, $\sum_{p\in\mathcal P}\norm{\Psi_{2,p}'f-\Psi_{2,p}'g}^2\le\norm{f-g}^2$, $\sum_{p}\norm{\Psi_{2,p}'f}^2\le\norm f^2$, and $\norm{\Phi^{(2)}f-\Phi^{(2)}g}\le\sqrt3\norm{f-g}$.
\end{proposition}
\begin{proof}
Since the $S_{\beta,j}$ have mutually orthogonal ranges (\cref{prop:struct}), $\abs{\abs a-\abs b}\le\abs{a-b}$ together with Minkowski yields $\sum_j\norm{\Psi_{1,j}'f-\Psi_{1,j}'g}^2\le\E_\beta\sum_j\norm{S_{\beta,j}(f-g)}^2=\E_\beta\norm{S^J_\beta(f-g)}^2\le\norm{f-g}^2$. Use this for $u_{j_1}=\Psi_{1,j_1}'f$ and $w_{j_1}=\Psi_{1,j_1}'g$ at the second layer to obtain $\sum_{j_2}\norm{\Psi_{1,j_2}'u_{j_1}-\Psi_{1,j_2}'w_{j_1}}^2\le\norm{u_{j_1}-w_{j_1}}^2$, noting that a sum restricted to a subset of $j_2$ is smaller. The second claim then results upon summing over $j_1$; setting $g=0$ in that claim yields the third. Finally, $\norm{A_J}\le1$ combined with Pythagoras across the three layers establishes the last claim.
\end{proof}

\begin{theorem}[Scale-separated channels]\label{thm:sep}
Take $\eps\le\frac12$ and set
\[
\mathcal E_{\rm sep}(f):=\Big(\sum_{j=0}^J\norm{\Psi_{1,j}'(L_\tau f)-L_\tau\Psi_{1,j}'f}^2\Big)^{1/2}.
\]
\begin{enumerate}[label=(\alph*),nosep]
\item \emph{Finite depth.} For Haar, $\mathcal E_{\rm sep}(f)\le\big(\sqrt2+(2+\sqrt{48})\sqrt{J+1}\big)\eps\norm f$; for $\psi\in C^1$, $\mathcal E_{\rm sep}(f)\le\big(\sqrt2+10L^2K_1\sqrt{J+1}\big)\eps\norm f$.
\item \emph{Affine warps.} If $\tau(x)=c+ax$ with $\abs a=\eps$, then $\mathcal E_{\rm sep}(f)\le C\eps\norm f$ for $\psi\in C^{1,\delta}$, whereas for Haar $\mathcal E_{\rm sep}(f)\le C\eps\big(1+\log_2\frac1\eps\big)^{1/2}\norm f$, uniformly in $J$.
\item \emph{$C^2$ warps.} If both $\psi\in C^{1,\delta}$ and $\tau\in C^2$ hold, then, uniformly in $J$,
\[
\mathcal E_{\rm sep}(f)\le C\big(\eps+\eps^{1-\delta}\curv^\delta+\eps^\delta\curv\big)\norm f .
\]
\end{enumerate}
\end{theorem}
\begin{proof}
\Cref{lem:local} holds at each scale; hence, with $Q_{\beta,j}f(x)=S_{\beta,j}[\tau'L_\tau f](x)+G_{\beta,j}f(x)$, we obtain $\abs{\Psi_{1,j}'(L_\tau f)(x)-L_\tau\Psi_{1,j}'f(x)}\le\E_\beta\abs{Q_{\beta,j}f(x)}$ as in \cref{prop:jac}. Applying Minkowski in $\ell^2(L^2(d\beta;L^2))$ together with orthogonality of the ranges of the $S_{\beta,j}$ yields
\[
\mathcal E_{\rm sep}(f)\le\big(\E_\beta\norm{S^J_\beta[\tau'L_\tau f]}^2\big)^{1/2}+\Big(\sum_j\E_\beta\norm{G_{\beta,j}f}^2\Big)^{1/2}\le\sqrt2\,\eps\norm f+\Big(\sum_j\E_\beta\norm{G_{\beta,j}f}^2\Big)^{1/2}.
\]
(a) For Haar, $G_{\beta,j}=M_{\beta,j}+R_{\beta,j}$ (the parts of \cref{prop:split} at scale $j$); for every $j$, the proof of \cref{prop:haar-MR} yields $\E_\beta\norm{M_{\beta,j}f}^2\le4\eps^2\norm f^2$ and $\E_\beta\norm{R_{\beta,j}f}^2\le48\eps^2\norm f^2$. Apply \cref{prop:C1} scale by scale for $\psi\in C^1$. Adding the $J+1$ equal bounds yields $\sqrt{J+1}$.

(b) $G_{\beta,j}=\tilde G_{\beta,j}\Lambda_\tau$; moreover, each $\tilde G_{\beta,j}$ is covariant as in \cref{lem:avg} when $\tau$ is affine. Taking $\bar\mu$ as in \eqref{eq:mubar}, we obtain $\E_\beta\norm{\tilde G_{\beta,j}g}^2\le\int\bar\mu(2^{-j}\xi)^2\abs{\hatt g(\xi)}^2d\xi$ from \cref{lem:avg}, together with $\sum_j\bar\mu(2^{-j}\xi)^2\le\sup\bar\mu\cdot\sum_j\bar\mu(2^{-j}\xi)\le C\eps^2$ for $\delta>0$. In the Haar case, at most $\log_2\frac1\eps+1$ terms of size $C\eps^2$ come from the range $1\le\abs s\le\eps^{-1}$. Since the remaining two ranges are geometric, the sum is $\le C\eps^2(1+\log_2\frac1\eps)$.

(c) Decompose $\tilde G_{\beta,j}=F_j-(R_j-R^\ell_j)-R^\ell_j$ as in \cref{ssec:C2}. By the proof of \cref{lem:freeze}, $\norm{R_j-R^\ell_j}$ is bounded by $C(\eps\curv^\delta2^{-j\delta}+\eps^\delta\curv2^{-j})$, a quantity square-summable in $j$. For $R^\ell_j=\int_{-\eps}^\eps D_{b,j}(\chi_b\,\cdot)\,db$, combining the symbol bound \eqref{eq:Dbsym} with \cref{lem:avg} at every scale yields $(\sum_j\E_\beta\norm{D_{b,j}h}^2)^{1/2}\le C\abs b^\delta\norm h$; applying Minkowski in $b$ then yields $C\eps^{1+\delta}\norm g$. For the first-order part, express $\tau(x)-\tau(v)=\tau'(x)(x-v)-r(x,v)$ with
\[
\abs{r(x,v)}\le\min\big(\tfrac12\curv\abs{x-v}^2,\,2\eps\abs{x-v}\big)\le2\curv^\delta\eps^{1-\delta}\abs{x-v}^{1+\delta}.
\]
Then $F_j=-\tau'\cdot E_j+\check R_j$, with $(x-v)\partial_vk_j(x,v)$ the kernel of $E_j$ and $r(x,v)\partial_vk_j(x,v)$ the kernel of $\check R_j$. On $\abs{x-v}\le L2^{-j}$ the latter is $\le C\curv^\delta\eps^{1-\delta}2^{j(1-\delta)}$; consequently, Schur's test yields $\norm{\check R_j}\le C\curv^\delta\eps^{1-\delta}2^{-j\delta}$, which is square-summable in $j$. Given $E_j$, denote the left endpoint of $I$ by $c_I$ and set $x-v=(x-c_I)-(v-c_I)$. As $\partial_vk_j(x,v)=\sum_I\Theta_I(x)\psi_I'(v)$,
\[
E_jg=\sum_{\abs I=2^{-j}}\Big[(s\Theta)_I\,\ip g{(\psi')_I}-\Theta_I\,\ip g{(t\psi')_I}\Big],
\]
where the $L^2$-normalized lattice copies of $s\Theta(s)$, $\psi'(t)$ and $t\psi'(t)$ are $(s\Theta)_I$, $(\psi')_I$ and $(t\psi')_I$, respectively. No more than $L+1$ supports overlap at any fixed scale; thus $\norm{E_jg}^2\le2(L+1)\big[\norm{s\Theta}^2\sum_{\abs I=2^{-j}}\abs{\ip g{(\psi')_I}}^2+\norm\Theta^2\sum_{\abs I=2^{-j}}\abs{\ip g{(t\psi')_I}}^2\big]$. The families $\{(\psi')_I\}$ and $\{(t\psi')_I\}$, taken over all $I\in\D^\beta$, are Bessel sequences [Gra14]: indeed, $\psi'$ and $t\psi'$ both have compact support, $C^\delta$ regularity and mean zero ($\int t\psi'=-\int\psi=0$). It follows that, for each $\beta$, $\sum_j\norm{E_jg}^2\le C\norm g^2$ and $\big(\sum_j\norm{F_jg}^2\big)^{1/2}\le C(\eps+\curv^\delta\eps^{1-\delta})\norm g$. Combining these, apply $\norm{\Lambda_\tau}\le\sqrt{3/2}$ and observe that $\eps\curv^\delta\le\eps^{1-\delta}\curv^\delta$ and $\eps^{1+\delta}\le\eps$.
\end{proof}

\begin{remark}
The distinction between the separated channels and the summed channel is twofold. Because scale-$j$ errors are measured in $\ell^2$ instead of being added, the loss at finite depth is only $\sqrt{J+1}$. A Bessel argument in (c) makes the first-order part square-summable on every lattice. This dispenses with the $T(1)$ theorem, at the cost of replacing $\eps\curv^\delta$ by the curvature term $\eps^{1-\delta}\curv^\delta$. For Haar with non-affine warps, the open question for \cref{conj:haar} also applies to the separated channels.
\end{remark}

\begin{corollary}[Depth two]\label{cor:depth2}
Suppose that, for every $u\in L^2$, $\mathcal E_{\rm sep}(u)\le B\norm u$, where $B$ is from any case of \cref{thm:sep}. It follows that
\[
\Big(\sum_{p\in\mathcal P}\norm{\Psi_{2,p}'(L_\tau f)-L_\tau\Psi_{2,p}'f}^2\Big)^{1/2}\le2B\norm f,\qquad\norm{\Phi^{(2)}(L_\tau f)-\Phi^{(2)}f}\le\sqrt3\,\big(2B+a_J(\tau)\big)\norm f .
\]
In particular, the depth-two map for $\psi\in C^{1,\delta}$ and $\tau\in C^2\cap L^\infty$ is stable, and its constant is independent of $J$.
\end{corollary}
\begin{proof}
Let $u_{j_1}=\Psi_{1,j_1}'f$ and $e_{j_1}=\Psi_{1,j_1}'(L_\tau f)-L_\tau u_{j_1}$; consequently, $\sum\norm{u_{j_1}}^2\le\norm f^2$ (\cref{prop:nonexp}) and $\sum\norm{e_{j_1}}^2\le B^2\norm f^2$. Pointwise,
\[
\abs{\Psi_{2,(j_1,j_2)}'(L_\tau f)-L_\tau\Psi_{2,(j_1,j_2)}'f}\le\E_{\beta'}\abs{S_{\beta',j_2}e_{j_1}}+\abs{\Psi_{1,j_2}'(L_\tau u_{j_1})-L_\tau\Psi_{1,j_2}'u_{j_1}} .
\]
Taking $\ell^2$ over $j_2$ yields an upper bound of $\norm{e_{j_1}}$ for the first term, by orthogonality, and of $B\norm{u_{j_1}}$ for the second. Next, taking $\ell^2$ over $j_1$ yields $2B\norm f$. For $\Phi^{(2)}$, the contribution from each layer consists of its covariance error together with the pooling error $a_J(\tau)$ multiplied by that layer's energy, namely $\le\norm f$; Pythagoras then combines the three layers.
\end{proof}
\noindent Induction shows that a depth-$M$ network constructed in this manner is non-expansive, with stability constant no greater than $MB$ before pooling, just as in wavelet scattering [M12].

\begin{corollary}[Dilation averaging]\label{cor:dil}
Let $\Pi$ denote an arbitrary probability measure on $[1,2)$. With $\curv$ replaced by $2\curv$, all covariance bounds established in \cref{sec:finiteJ,sec:uniform} for $\Psi_1'$ and in \cref{thm:sep} for the scale-separated channels remain valid for the dilation-averaged channels, with no change in the constant. Neither $\eps$ nor the affine class changes. Non-expansiveness and the pooling bound likewise remain unaffected.
\end{corollary}
\begin{proof}
Set $\tau_r(x)=\tau(rx)/r$. Computing directly yields $D_r^{-1}L_\tau=L_{\tau_r}D_r^{-1}$ as well as $\abs{D_ru}=D_r\abs u$. Consequently, $\E_\beta\abs{S^J_{\beta,r}\,\cdot\,}=D_r\Psi_1'D_r^{-1}$ and
\[
\E_\beta\abs{S^J_{\beta,r}L_\tau f}-L_\tau\E_\beta\abs{S^J_{\beta,r}f}=D_r\big[\Psi_1'L_{\tau_r}-L_{\tau_r}\Psi_1'\big]D_r^{-1}f .
\]
$D_r$ is unitary; moreover, $\norm{\tau_r'}_\infty=\eps$, $\norm{\tau_r''}_\infty=r\curv\le2\curv$, $\norm{\tau_r}_\infty\le\norm\tau_\infty$, and whenever $\tau$ is affine, so is $\tau_r$. Thus the bound for $\tau_r$ holds for every fixed $r$. Now integrate over $r$, using $\abs{\int a\,d\Pi-\int b\,d\Pi}\le\int\abs{a-b}\,d\Pi$ together with Minkowski's inequality. For the scale-separated channels, apply this reasoning to each channel separately.
\end{proof}

\subsection{Haar with general warps}
For Haar, one cannot use \cref{prop:F,lem:Rj}, since the Haar kernel on a single lattice is not Calder\'on--Zygmund. Nor does the superposition of \cref{lem:super} work: because $\Gamma_y(s)$ has no decay for Haar, the symbol of every $D_b$ is non-summable. The $b$-integration resolves this in the affine case by reassembling the bracket $F(s)-\lambda F(\lambda s)$. That cancellation is broken by the cutoffs $\chi_b(v)$ when $\tau'$ varies.

\begin{conjecture}[Haar, general warps]\label{conj:haar}
For Haar and $\tau\in C^2$, the bound $(\E_\beta\norm{G_\beta f}^2)^{1/2}\le C(\eps\log\frac e\eps+\curv)\norm f$ holds uniformly in $J$. It may be possible to omit the logarithm.
\end{conjecture}

\noindent One natural approach freezes the slope within windows of size $w$ and uses \cref{thm:B} separately in each, with the averaged strip bounds of \cref{prop:haar-MR} controlling the slope variation $\curv w$. What remains missing is almost-orthogonality across scales for the errors after averaging. Numerical evidence comes from \Cref{ssec:num-worst}: for both warps, the averaged worst case for Haar saturates in $J$. As $\eps$ drops from $0.1$ to $0.02$, this quantity stays unchanged for a nine-scale warp but increases by approximately $10\%$ for a smooth warp.

%==================================================================
\section{Numerical experiments}\label{sec:num}

\subsection{Implementation}\label{ssec:impl}
We use PyWavelets' periodized transform on a periodic grid of $N_{\rm s}$ samples; the discrete child offset $m_d$ corresponds to the continuous offset $m=m_d+(L-1)/2$. Shifts by $2^{L_{\max}}$ samples leave the lattice invariant, so the lattice expectation is computed \emph{exactly} by averaging over the $2^{L_{\max}}$ integer offsets (cycle spinning [CD95]). All constants are taken from the infinite-product spectrum rather than from sampled wavelets, whose values drift with the sampling level. The reference implementation passes 31 regression checks, and \Cref{tab:aligned} and the certified ripples were recomputed with independent code. Further details (offset matching, dilation quadrature, sources of numerical error and the regression checks) are given in Section~S1 of the supplementary material, and the scripts are listed in \cref{app:repro}. As stated in the declaration at the end of the paper, generative AI tools were used to debug the code.

\subsection{Deformation experiments: design}\label{ssec:design}
The number of samples in the finest target cell is fixed (16 or 8), with $J+1$ shift scales above that cell. Scale invariance makes this equivalent to refinement with the coarsest cell held fixed. All $2^{L_{\max}}$ integer offsets, $L_{\max}=J+4$ or $J+5$, are included in the exhaustive lattice average. Because the inputs are trigonometric polynomials, evaluation of $L_\tau f$ achieves relative accuracy $10^{-10}$ using $16\times$ FFT oversampling followed by a periodic cubic spline. Two warps are used: \emph{smooth}, $\tau=\eps\frac P{2\pi}\sin\frac{2\pi x}P$ with $P=N_{\rm s}/8$ or $N_{\rm s}/4$; and \emph{multi-scale}, which sums eight or nine sines having random phases and periods of $2^4,\dots,2^{11}$ or $2^{12}$ samples. For each component, the derivative amplitude is $\eps$ divided by the number of components, giving $\norm{\tau'}_\infty\le\eps$. Large curvature at the finest scales makes the multi-scale warp a probe of the regime of \cref{conj:haar}. Haar and aligned sym4 serve as the generators. Measuring the discretization floor with a constant fractional shift $\tau\equiv0.5$, for which zero is predicted by the theory, gives at most $3.5\times10^{-4}$.

\subsection{Typical inputs}\label{ssec:num-typ}
\begin{figure}[tbp]\centering
\includegraphics[width=0.93\linewidth]{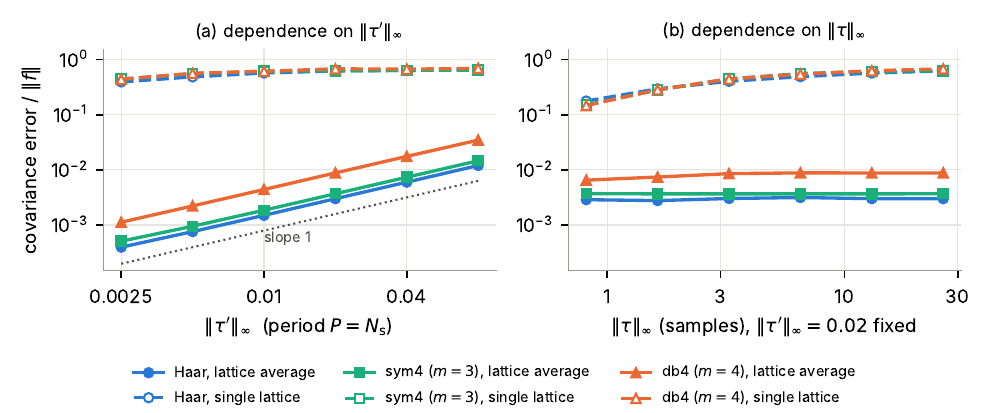}
\caption{Covariance error with $2^{13}$ samples, $J=5$, exhaustive average over 512 offsets and $\tau=\eps\frac P{2\pi}\sin\frac{2\pi x}P$. Generators: Haar and the aligned sym4 ($m=3$, discrete offset $m_d=0$) and db4 ($m=4$, $m_d=1$). (a) The lattice average is linear in $\eps$ for all three, while a single lattice is $O(1)$. (b) Increasing $\norm\tau_\infty$ by a factor $32$ at fixed $\norm{\tau'}_\infty$ changes the averaged error by at most $2\%$ (sym4), $14\%$ (Haar) and $35\%$ (db4, constant from $\norm\tau_\infty\approx6$ samples on), while the single-lattice error grows by a factor $3.5$--$4.6$.}\label{fig:deform}
\end{figure}

\Cref{tab:mod} demonstrates that, for both generators and both warps across every input class, the covariance error remains flat in $J$. These runs also establish two points relevant to the proofs. (i) When all seven scales are computed together rather than each scale separately, the total is \emph{smaller} than even $(\sum_je_j^2)^{1/2}$, the orthogonal sum of the single-scale errors. Scales combine sub-additively through the modulus, without coherent accumulation. (ii) After normalization by the energy of $f$ within the active scale window, the \emph{linear} defect $(\E_\beta\norm{Q_\beta f}^2)^{1/2}$ of \eqref{eq:reduce2} remains constant in $J$: for Haar it is $0.12\pm0.003$, and for sym4 it is $0.13$--$0.17$ at $\eps=0.1$. This linear defect is what the proofs of \cref{sec:uniform} use; they do not take advantage of the additional cancellation within the modulus, amounting to a further factor of $3$--$7$. For $J=5$, \Cref{fig:deform} demonstrates linearity of the averaged error in $\eps$, together with boundedness when $\norm\tau_\infty$ increases $32$-fold while $\norm{\tau'}_\infty$ is fixed, in agreement with the prediction of \cref{thm:cov}. For the aligned db4, the constant is about $2.4$ times the sym4 constant, whereas the linear dependence on $\eps$ is unchanged; both generators are $C^{1,\delta}$.

\begin{table}[htbp]\centering\footnotesize\setlength{\tabcolsep}{3.2pt}
\resizebox{\linewidth}{!}{\begin{tabular}{lll*{7}{c}}\toprule
gen. & warp & input & $J=0$ & 1 & 2 & 3 & 4 & 5 & 6\\\midrule
Haar & smooth & white & $.0298\pm.0008$ & $.0308\pm.0011$ & $.0309\pm.0010$ & $.0310\pm.0009$ & $.0309\pm.0008$ & $.0309\pm.0008$ & $.0308\pm.0008$\\
Haar & smooth & pink & $.0154\pm.0009$ & $.0161\pm.0010$ & $.0162\pm.0009$ & $.0165\pm.0009$ & $.0163\pm.0009$ & $.0160\pm.0007$ & $.0154\pm.0009$\\
Haar & multi & white & $.0088\pm.0003$ & $.0092\pm.0003$ & $.0090\pm.0003$ & $.0091\pm.0003$ & $.0091\pm.0003$ & $.0090\pm.0003$ & $.0090\pm.0003$\\
Haar & smooth & balanced & .0345 & .0344 & .0339 & .0307 & .0266 & .0241 & .0242\\
sym4 & smooth & white & $.0373\pm.0018$ & $.0395\pm.0014$ & $.0390\pm.0015$ & $.0386\pm.0012$ & $.0383\pm.0012$ & $.0383\pm.0011$ & $.0381\pm.0012$\\
sym4 & smooth & pink & $.0214\pm.0015$ & $.0228\pm.0013$ & $.0219\pm.0015$ & $.0223\pm.0012$ & $.0210\pm.0013$ & $.0201\pm.0009$ & $.0192\pm.0009$\\
sym4 & multi & white & $.0110\pm.0003$ & $.0119\pm.0002$ & $.0117\pm.0002$ & $.0115\pm.0003$ & $.0114\pm.0003$ & $.0113\pm.0003$ & $.0113\pm.0003$\\
sym4 & smooth & balanced & .0415 & .0467 & .0396 & .0349 & .0300 & .0293 & .0284\\\midrule
Haar & \multicolumn{2}{l}{floor ($\tau\equiv0.5$)} & .0003 & .0003 & .0002 & .0001 & .0001 & .0001 & .0000\\
sym4 & \multicolumn{2}{l}{floor ($\tau\equiv0.5$)} & .0005 & .0004 & .0003 & .0002 & .0001 & .0001 & .0001\\\bottomrule
\end{tabular}}
\caption{Covariance error $\norm{\Psi_1'(L_\tau f)-L_\tau\Psi_1'f}/\norm f$ against $J$ for $\eps=0.1$ ($N_{\rm s}=2^{14}$, finest cell 16 samples, mean $\pm$ s.d.\ over 3 seeds). Inputs: white (flat spectrum), pink (equal energy per octave), balanced (one wavelet per active level at a common point). Every row is flat in $J$, and doubling $\eps$ doubles the error.}\label{tab:mod}
\end{table}

\subsection{Scale-separated channels, a second layer, and dilation averaging}\label{ssec:num-ext}
We repeated the typical-input experiment for the variants of \cref{ssec:ext} (smooth warp, $\eps=0.1$, one white and one pink input, $J=0,\dots,6$), measuring the error per unit of input energy in the active scale window; the table and further details are in Section~S2 of the supplementary material. The scale-separated error $\rho_{\rm sep}$ stays flat in $J$ (within $2\%$ for white inputs, slowly decreasing for pink inputs), so these inputs do not reach the factor $\sqrt{J+1}$ of \cref{thm:sep}(a). The second-layer error is $1.4$--$1.9$ times $\rho_{\rm sep}$, within the factor $2$ of \cref{cor:depth2}, and the second-layer energy is $1$--$8\%$ of the first-layer energy. With four dilation nodes and the uniform measure, each node stays within $25\%$ of the $r=1$ value, as \cref{cor:dil} predicts, and the dilation-averaged error is $60$--$78\%$ of the translation-averaged error, again flat in $J$.

\subsection{Worst-case inputs}\label{ssec:num-worst}
\paragraph{Method.} We determine the largest singular value $\sigma_{\max}$ of the linear map $A:f\mapsto2^{-L_{\max}/2}(Q_bf)_b$ restricted to the space of band-limited inputs (periods $\ge16$ samples) and allowing $b$ to range over all lattice offsets. It follows from \eqref{eq:reduce2} that $\sigma_{\max}$ is the best constant for the linear bound $(\E_b\norm{Q_bf}^2)^{1/2}\le\sigma_{\max}\norm f$ on band-limited inputs; this is precisely the quantity controlled by the proofs of \cref{sec:uniform}. Our computation uses $N_{\rm s}=2^{13}$ and finest target cells of 8 samples, together with the exact local translation identity $S_bT_{n+\theta}f(x)=S_{b-n}(T_\theta f)(x-n)$, exact FFT translates on seven $\theta$-nodes, and 4-point Lagrange interpolation. We implement the adjoint $A^*$ explicitly. The relation $\ip{Af}Z=\ip f{A^*Z}$ holds to $2\times10^{-15}$, while a constant shift yields $Af=0$ to machine precision. To determine the largest eigenvalue of $A^*A$, we apply the Lanczos method with full reorthogonalization, starting from a scale-balanced input (one wavelet per active level at a common point) combined with $5\%$ band-limited Gaussian noise with seed $J$. After a minimum of four steps, iteration terminates once the residual of the top Ritz pair is below $10^{-9}\theta$, or upon reaching $40$ steps. We retain $\theta^{1/2}$ and the residual $r=\norm{A^*Ay-\theta y}$ for the top Ritz pair $(\theta,y)$. The Rayleigh quotient property of $\theta$ implies $\theta\le\lambda_{\max}$, with some eigenvalue belonging to $[\theta-r,\theta+r]$. The relative residual $r/\theta$ never exceeds $1\times10^{-3}$ and is below $10^{-6}$ for most entries. A larger eigenvalue with an eigenvector almost orthogonal to the starting vector cannot, however, be ruled out by Lanczos alone. For this reason, we additionally constructed the full $1024\times1024$ matrix of $A^*A$ on the band-limited space and calculated all its eigenvalues (Haar for $J\le3$ and sym4 for $J\le2$, each using the smooth warp and $\eps=0.1$; symmetry of the matrices holds to $6\times10^{-16}$). Across all cases, the largest eigenvalue matches the Lanczos value to a relative $8\times10^{-6}$ (Table~S2 in the supplementary material). Entries with $J=8$ instead come from power-iteration values (10--14 steps) obtained in an earlier run. Being Rayleigh quotients, these values are lower bounds; whenever both results are available, they fall $0.1$--$1.5\%$ below the Lanczos values. Every reported value concerns the discrete operator, with band-limited inputs and an exhaustive average over integer lattice offsets.

\begin{table}[h]\centering\small\setlength{\tabcolsep}{4pt}
\begin{tabular}{lllccccccccc}\toprule
gen. & warp & $\eps$ & $J=0$ & 1 & 2 & 3 & 4 & 5 & 6 & 7 & 8\\\midrule
Haar & smooth & 0.1 & 1.294 & 1.508 & 1.945 & 1.977 & 1.995 & 2.008 & 2.017 & 2.018 & ---\\
Haar & smooth & 0.3 & 1.367 & 1.766 & 1.820 & 1.938 & 1.953 & 1.962 & 1.963 & --- & ---\\
Haar & smooth & 0.02 & 1.255 & --- & --- & --- & 2.137 & --- & --- & --- & 2.208$^\ast$\\
\addlinespace[2pt]
sym4 & smooth & 0.1 & 2.073 & 2.388 & 2.386 & 2.386 & 2.387 & 2.387 & 2.387 & --- & ---\\
sym4 & smooth & 0.3 & 2.110 & 2.776 & 2.740 & 2.749 & 2.749 & 2.749 & 2.749 & --- & ---\\
\addlinespace[2pt]
Haar & multi & 0.1 & 0.806 & --- & 1.225 & --- & 1.313 & --- & 1.338 & --- & 1.323$^\ast$\\
Haar & multi & 0.05 & 0.800 & --- & 1.204 & --- & 1.302 & --- & 1.335 & --- & 1.326$^\ast$\\
Haar & multi & 0.02 & 0.796 & --- & 1.189 & --- & 1.291 & --- & 1.327 & --- & 1.319$^\ast$\\
\bottomrule\end{tabular}
\caption{Worst-case linear defect $\sigma_{\max}/\norm{\tau^\prime}_\infty$ of the lattice-averaged shift ($2^{13}$ samples, finest cell 8 samples). Unmarked entries are Lanczos values, with relative residual at most $1\times10^{-3}$ (\cref{ssec:num-worst}); entries marked ${}^\ast$ are power-iteration lower bounds. For the multi-scale warp the normalization uses $\norm{\tau^\prime}_\infty=0.747\,\eps$, the maximum of $\abs{\tau^\prime}$ over the $2^{13}$ sample points; the supremum over $\R$ is $0.749\,\eps$, so the normalized values change by less than $0.4\%$. Every row saturates in $J$; the three multi-scale Haar rows agree to within $3\%$ at every $J$.}\label{tab:worst}
\end{table}

\begin{figure}[tbp]\centering
\includegraphics[width=\linewidth]{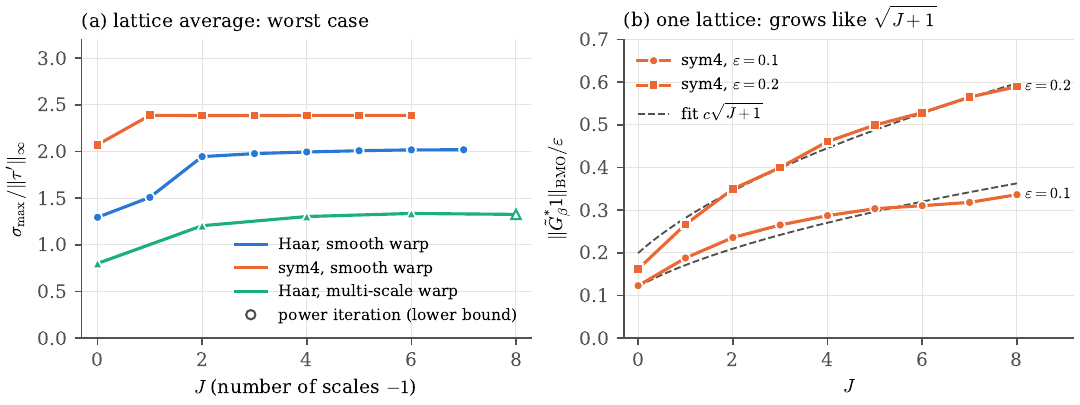}
\caption{Lattice average versus a single lattice. (a) Averaged worst-case linear defect $\sigma_{\max}/\norm{\tau'}_\infty$ (Lanczos values, filled markers; power-iteration lower bounds, open markers; \cref{ssec:num-worst}): flat in $J$ for Haar and aligned sym4 with a smooth warp ($\eps=0.1$), and for Haar with the multi-scale warp ($\eps=0.05$). (b) Single-lattice $\norm{\tilde G^*_\beta1}_\BMO/\eps$ for aligned sym4 with a smooth warp, which keeps growing; the dashed curves are least-squares fits $c\sqrt{J+1}$.}\label{fig:uniform}
\end{figure}

\paragraph{Findings.} (i) After only a few scales, $\sigma_{\max}/\norm{\tau'}_\infty$ saturates for both generators and both warps (\cref{tab:worst} and \cref{fig:uniform}(a)). Its increments decrease geometrically; starting at $J=3$, the change in any row is at most about $2\%$ (less than $0.1\%$ for sym4). The smooth warp at $\eps=0.02$ is an exception, as discussed in (iv). Neither $(J+1)$ nor $\sqrt{J+1}$ growth is apparent. (ii) Only a factor of $1.4$--$1.8$ separates the worst case from the typical linear defect. (iii) Rather than point concentrations, the maximizers are fine-scale oscillations: they extend across hundreds of samples and have spectral centroid at the finest active scales. (iv) Dependence on $\eps$ for Haar. For every $J$, the multi-scale warp gives normalized values at $\eps=0.02$, $0.05$ and $0.1$ that differ by no more than $3\%$. For the smooth warp, reducing $\eps$ from $0.1$ to $0.02$ increases the value at $J=4$ from $1.99$ to $2.14$ (Lanczos). The largest computed value also increases, from $2.02$ ($J=7$, Lanczos) to at least $2.21$ ($J=8$, power-iteration lower bound); saturation then occurs later, at $J\approx\log_2\frac1\eps+2$. At this point the range $1\le\abs s\le\eps^{-1}$ of \eqref{eq:mubar}, which is the source of the logarithm in \cref{thm:B}, becomes fully populated. An entire factor $\log\frac e\eps$ would increase the value by $49\%$. Thus, although the data support $J$-uniformity in \cref{conj:haar}, they leave the presence of a weak logarithm unresolved.

\subsection{Single lattice versus lattice average}\label{ssec:num-single}
On one lattice, direct evaluation gives
\[
h_\beta=\sum_{j\le J}h_j,\qquad h_j(u)=\sum_{\abs I=2^{-j}}\int\Theta_I(x)\,\psi_I(u+\tau(x))\,dx,
\]
which is equal to $-\tilde G^*_\beta1\circ\rho^{-1}$, the function with BMO norm appearing in the $T(1)$ argument. Our periodic domain consists of four unit cells with lattice offset $\beta=0.3137$. For $\Theta$ we employ a quadrature with weights that reproduce its vanishing moments exactly; the BMO norm is taken to be the $L^2$ mean oscillation across all dyadic windows. For the averaged setting, we obtain the exact constants $\sup_\xi\M_J(\xi)^{1/2}/\eps$ for affine warps using \cref{lem:avg,lem:sym}, with exact $\hatt\psi$, 32 samples per finest cell and $\xi\in[2^{-2},2^{16}]$.

\begin{table}[htbp]\centering\small
\begin{tabular}{lll*{5}{c}}\toprule
 & gen. & $\eps$ & $J=0$ & 2 & 4 & 6 & 8\\\midrule
single lattice: & sym4 & 0.05 & 0.093 & 0.182 & 0.235 & 0.256 & 0.256\\
$\norm{h_\beta}_\BMO/\eps$ & sym4 & 0.1 & 0.123 & 0.236 & 0.287 & 0.310 & 0.336\\
(smooth warp) & sym4 & 0.2 & 0.162 & 0.349 & 0.461 & 0.529 & 0.590\\
 & db4 & 0.1 & 0.084 & 0.180 & 0.231 & 0.257 & 0.278\\
 & Haar & 0.1 & 3.214 & 4.112 & 4.220 & 4.270 & 4.766\\\bottomrule
\end{tabular}
\caption{Single-lattice BMO norm of $\tilde G_\beta^*1$ (up to the change of variables $\rho$). For sym4 and db4 it grows like $\sqrt{J+1}$, with an amplitude of order $\eps^{1.3}$--$\eps^{1.6}$. For Haar, $\norm{h_\beta}_\BMO\approx0.3$--$0.5$ at every $\eps$ tested ($0.05$--$0.2$), i.e.\ it is $O(1)$.}\label{tab:bmo}
\end{table}

\begin{table}[htbp]\centering\small
\begin{tabular}{llcccc}\toprule
generator & $\eps=\abs a$ & $J=2$ & $J=5$ & $J=8$ & $J=11$\\\midrule
Haar & 0.05 & 1.8598 & 2.3879 & 2.4928 & 2.4945\\
Haar & 0.1 & 1.9217 & 2.3252 & 2.3665 & 2.3750\\
Haar & 0.3 & 2.2811 & 2.4769 & 2.5101 & 2.5154\\
sym4 & 0.05 & 2.5809 & 2.5816 & 2.5816 & 2.5816\\
sym4 & 0.1 & 2.6930 & 2.6936 & 2.6936 & 2.6936\\
sym4 & 0.3 & 3.2213 & 3.2223 & 3.2223 & 3.2223\\\bottomrule
\end{tabular}
\caption{Exact lattice-averaged constants $\sup_\xi\M_J(\xi)^{1/2}/\eps$ for affine warps $\tau=ax$ (\cref{thm:B}). sym4 is constant in $J$ to four digits. Haar saturates by $J\approx8$, and its constant hardly depends on $\eps$.}\label{tab:aff}
\end{table}

\Cref{tab:bmo,tab:aff,fig:uniform} exhibit the dichotomy of \cref{sec:uniform}. The remainder accumulates when a single lattice is used. With the affine warp $\tau=\eps x$ and aligned sym4, dividing the root mean square of $h_\beta$ by $\eps^{1.4}\sqrt{J+1}$ gives $0.165\pm0.002$ for $J=0,\dots,10$; this agrees with the exact lacunary formula of \cref{thm:single-lower} (\cref{rem:sharp}). Averaging over lattices, by contrast, leaves the constants flat. Under a multi-scale warp, growth on a single lattice is first order in $\eps$ until every component of $\tau'$ has been resolved, at which point it saturates, in keeping with the Carleson structure of $F^*1$.

%==================================================================
\section{Discussion}\label{sec:disc}

\paragraph{Summary.}
For the averaged DSSN channel, the local translation identity absorbs translations exactly, while the Jacobian incurs a $J$-free $\sqrt2\eps$; what remains of the defect is a dyadic commutator. With $C^{1,\delta}$ generators, the first-order part of this commutator is uniformly bounded on every lattice. Although its second-order part accumulates on a single lattice whenever $\sigma_\eps>0$, averaging over lattices turns it into a decaying Fourier multiplier. Consequently, stability is $J$-uniform both for affine warps (Haar and $C^{1,\delta}$ generators, with a logarithm for Haar) and for $C^2$ warps ($C^{1,\delta}$ generators). The case of Haar with general warps is still open, although the computations of \cref{ssec:num-worst} agree with the conjectured bound. These conclusions also apply to scale-separated channels (using $\sqrt{J+1}$ in place of $J+1$ at finite depth), to a second layer, and to dilation-averaged lattices.

\paragraph{Limitations.} (1) \Cref{thm:single-lower} supplies a lower bound for the single-lattice defect given by a multiple of $\sigma_\eps$. For sym4 and db4, positivity of $\sigma_\eps$ has been established numerically only, at three values of $\eps$. (2) For Haar, a proof of the $J$-uniform bound is available only for affine warps. With general warps, it is \cref{conj:haar}; the computations for $J\le8$ support it but do not prove it. (3) The $T(1)$ theorem in \cref{sec:uniform} yields constants that are not explicit, while the curvature terms in \cref{thm:C} are possibly nonoptimal. (4) Every result is one-dimensional and pertains to the complexity-one shift. (5) The treatment of the network is analytic, with no claim regarding learning or classification performance.

\paragraph{Further directions.}
The principal question that remains open, \cref{conj:haar}, concerns a $J$-uniform bound for Haar with general warps and the necessity of the logarithm. Another question concerns removing the curvature terms from \cref{thm:C} and determining whether lattice averaging alters the $C^{1+}$ threshold identified by Nicola--Trapasso [NT23] for wavelet scattering. Although our preliminary runs for finite ensembles of $M$ lattices agree with an $M^{-1/2}$ fluctuation in norm, a proof is lacking. An explanation of the aligned offset $m^\ast$ should come from the group delay of the filter, along the lines of Selesnick's condition [S01]. For general Calder\'on--Zygmund operators, the shifts of higher complexity in [HL22] are required. A natural next step is therefore to extend the stability theory to channels constructed from these higher-complexity shifts, including the decay of second-layer energy with path length.

\section*{Declaration of generative AI and AI-assisted technologies in the manuscript preparation process}
During the preparation of this work the author used Claude (Anthropic; mainly Claude Opus~5.5) and ChatGPT (OpenAI; mainly GPT-6) in three ways: in the early stage, for preliminary numerical experiments that helped to shape the concept; during the project, to verify and debug the experiment code, whose results the author then checked again; and to improve the language of the manuscript. After using these tools, the author reviewed and edited the content as needed and takes full responsibility for the content of the published article.

\section*{Acknowledgments}
This work was supported by the National Research Foundation of Korea (NRF) grant funded by the Korea government (MSIT) (RS-2026-25602221).

%==================================================================
\appendix
\section{Reproducibility}\label{app:repro}
All computations were done in Python with NumPy, SciPy and PyWavelets, in double precision unless stated otherwise. The most recent computations (\cref{tab:worst,fig:deform}, Tables~S1 and~S2 of the supplementary material, the certified ripples and the recomputation of \cref{tab:aligned}) used Python 3.13.18, NumPy 2.5.3, SciPy 1.18.1, PyWavelets 1.8.0, mpmath 1.4.1 (interval arithmetic, 30 digits) and Matplotlib 3.11.2. The earlier computations were reproduced with this environment where they were rerun: the Haar curves of \cref{fig:deform}, the $J=0$ entries of \cref{tab:mod} and all entries of \cref{tab:aligned} agree to the digits shown. All random inputs use fixed seeds, stated in the scripts. The scripts and their outputs are available at \url{https://github.com/chdaewon/random-dyadic-shift-scattering} and are archived at Zenodo [Chu26]; the paths below are relative to the repository root.
\begin{center}\small
\begin{tabular}{>{\raggedright\arraybackslash}p{5.2cm}>{\raggedright\arraybackslash}p{9.6cm}}\toprule
result & script and output\\\midrule
\cref{cor:transfail}, ripples, \cref{tab:aligned} & \texttt{verify/verify\_constants.py}, \texttt{verify\_constants.json}\\
\cref{tab:mod}, single-scale and linear-defect checks & \texttt{theory\_J/ej\_check.py}, \texttt{ej\_run.py}, \texttt{ej\_linear.py}, \texttt{ej\_run.json}, \texttt{ej\_linear.json}\\
\cref{tab:worst}, \cref{fig:uniform}(a) & \texttt{theory\_J/lz\_worst.py} \texttt{<gen> <eps> <warp> <Js> <steps> <full-Js>} (Lanczos and full matrix), \texttt{lz\_*.json}; power iteration: \texttt{theory\_J/ej\_worst.py}, \texttt{ej\_worst\_*.json}; table: \texttt{paper\_final/make\_tab\_worst.py}\\
\cref{cor:transfail}, ripples of \cref{tab:ripple} (Haar; certified) & \texttt{paper\_final/review/cert\_ripple.py}, \texttt{cert\_ripple.json}\\
\cref{tab:aligned} (independent recomputation); sign check of \cref{prop:mult}(c) & \texttt{paper\_final/review/recompute\_tab1.py}, \texttt{check\_sign.py}\\
Table~S2 (supplementary) & \texttt{theory\_J/lz\_worst.py} (option \texttt{<full-Js>}), \texttt{lz\_haar\_smooth\_eps0.1.json}, \texttt{lz\_sym4\_smooth\_eps0.1.json}\\
Table~S1 (supplementary) & \texttt{theory\_J/ej\_ext.py} \texttt{layers}\textbar\texttt{dil <gen> <Js>}, \texttt{ej\_ext\_*.json}; table: \texttt{paper\_final/make\_tab\_ext.py}\\
\cref{tab:bmo}, \cref{fig:uniform}(b) & \texttt{theory\_J/ej\_bmo.py}; \texttt{ej\_bmo\_*.json}\\
\cref{rem:sharp} ($\sigma_\eps$, exact $\norm{S_J}$) & \texttt{theory\_J/ej\_single\_lowerbound.py}; \texttt{ej\_single\_lowerbound\_*.json}\\
\cref{tab:aff} & \texttt{theory\_J/ej\_affine\_symbol.py}; \texttt{ej\_affine\_symbol\_*.json}\\
low-frequency control runs & \texttt{theory\_J/ej\_lowfreq2.py}; \texttt{ej\_lowfreq2\_*.json}\\
\cref{fig:uniform} & \texttt{paper\_final/make\_fig5.py}\\
\cref{fig:deform} & \texttt{paper\_final/fig3/e3\_deform\_aligned.py}, \texttt{e3\_aligned.json}, \texttt{make\_fig3.py}\\
\cref{fig:ripple,fig:sample}, implementation checks & \texttt{early\_scripts/}; \texttt{regression/} (\texttt{dssn\_core} 0.2.0, 31 regression checks)\\\bottomrule
\end{tabular}
\end{center}

\clearpage
\appendix\renewcommand{\thesection}{S\arabic{section}}\setcounter{section}{0}
\setcounter{table}{0}\renewcommand{\thetable}{S\arabic{table}}
\renewcommand{\theHsection}{S\arabic{section}}\renewcommand{\theHtable}{S\arabic{table}}
\begin{center}\Large\textbf{Supplementary material}\end{center}
\noindent This supplement collects implementation details and two tables of numerical results that support \cref{sec:num}. It is submitted to the journal as a separate file; here it is appended to the preprint.

\section{Implementation details}\label{supp:impl}
\paragraph{Discrete shift.} We employ PyWavelets' periodized transform on a periodic grid containing $N_{\rm s}$ samples. Under the discrete shift, the detail coefficient $(j,k)$ maps to $(j+1,2k+m_d)$ and $(j+1,2k+m_d+1)$ with weights $\pm1/\sqrt2$. Because periodization re-centres the atoms, the continuous offset $m=m_d+(L-1)/2$ corresponds to the discrete offset $m_d$. To verify this correspondence, we matched translation-averaged multipliers for db2, db4, sym4, sym6 and coif3, obtaining a match error no greater than $8\times10^{-3}$; the runner-up offset performs worse by at least $50\times$. Thus, for the near-linear-phase families, $m_d=0$ is already aligned, as is $m_d^\ast=1$ for db4.

\paragraph{Exhaustive lattice averaging.} The coarsest cell contains $2^{L_{\max}}$ samples when there are $L_{\max}$ levels; shifts by $2^{L_{\max}}$ samples leave the lattice invariant. Consequently, the discrete lattice expectation \emph{exactly} is computed by averaging over the $2^{L_{\max}}$ integer offsets (cycle spinning [CD95]). Using $N_{\rm s}=4096$ and $64$ offsets, this computation is approximately $15\times$ faster than Monte Carlo using $1024$ draws.

\paragraph{Dilation quadrature.} For $Q=4$ and $N_{\rm s}=2^{14}$, the operator $A_Qf=\frac1Q\sum_qR_{r_q}^{-1}\E_\beta S_\beta R_{r_q}f$ reproduces the symbol of aligned sym4 to $7.5\times10^{-6}$ in the real part and $4.5\times10^{-5}$ in the imaginary part; here $R_r$ is an FFT resampling to a multiple of $2^{L_{\max}}$ samples.

\paragraph{Sources of numerical error.} Three effects may appear to disagree with the theory. With Haar, endpoint sampling introduces a spurious identity component of $0.12\%$; midpoint sampling eliminates it. Nodes are misplaced when the quadrature grid is not divisible by $Q$. Most significant is the sampling-level drift in constants obtained from sampled wavelets (\texttt{wavefun}): sym4 gives $\operatorname{Re}\kappa=-0.0074,-0.0093,-0.0103,-0.0108$ at levels $10$--$13$ and a Richardson limit of $-0.01130$. The latter agrees with the infinite-product value reported in \cref{tab:aligned}, the value used throughout.

\paragraph{Code verification.} The reference implementation undergoes 31 regression checks covering the discrete shift's adjointness, $S1=0$, the Haar constant $\kappa_h$ together with its ripples for $Q=1,2,4,5$, aligned constants for db4 and sym4, unaligned constants $\abs{c_0}$, signal-level dilation quadrature, and the discrete-to-continuous offset map. A single check is expected to fail because its comparison uses the value at sampling level $10$, thereby recording the sym4 drift of $\operatorname{Re}\kappa$ discussed above. Independent code was used to recompute \Cref{tab:aligned} and the certified ripples. Scripts and outputs are listed in \cref{app:repro}.

\section{Scale-separated channels, a second layer, and dilation averaging}\label{supp:ext}
We carry out the typical-input experiment from \cref{tab:mod} again, now using the variants of \cref{ssec:ext} with the smooth warp, $\eps=0.1$, one white input and one pink input, and $J=0,\dots,6$. Here, normalization by $\norm f$ gives a misleading picture: increasing $J$ brings more input energy into the active scale window, so a $J$-uniform bound $\mathcal E\le B\norm f$ permits growth of the error until all the energy lies within that window. Accordingly, the errors we report are per unit of energy within the relevant window (\cref{tab:ext}).

\begin{table}[t]\centering\small
\begin{tabular}{lll ccccccc}\toprule
gen. & input & & $J=0$ & 1 & 2 & 3 & 4 & 5 & 6\\\midrule
Haar & white & $\rho_{\rm sep}$ & .063 & .062 & .062 & .062 & .062 & .062 & .062\\
 &  & $\rho_2$ & --- & .104 & .095 & .098 & .096 & .102 & .106\\
 &  & $e_\Pi/e$ & .712 & .720 & .723 & .720 & .722 & .723 & .721\\
Haar & pink & $\rho_{\rm sep}$ & .047 & .046 & .044 & .043 & .044 & .045 & .043\\
 &  & $\rho_2$ & --- & .087 & .078 & .076 & .070 & .071 & .073\\
 &  & $e_\Pi/e$ & .760 & .777 & .718 & .686 & .667 & .634 & .657\\\midrule
sym4 & white & $\rho_{\rm sep}$ & .075 & .074 & .074 & .073 & .073 & .073 & .073\\
 &  & $\rho_2$ & --- & .106 & .122 & .121 & .117 & .121 & .126\\
 &  & $e_\Pi/e$ & .654 & .693 & .687 & .682 & .681 & .679 & .678\\
sym4 & pink & $\rho_{\rm sep}$ & .065 & .063 & .062 & .061 & .058 & .056 & .054\\
 &  & $\rho_2$ & --- & .102 & .119 & .114 & .111 & .107 & .104\\
 &  & $e_\Pi/e$ & .698 & .717 & .707 & .625 & .625 & .601 & .622\\\bottomrule
\end{tabular}
\caption{Scale-separated channels, second layer and dilation averaging ($N_{\rm s}=2^{14}$, finest cell 16 samples, smooth warp, $\eps=0.1$). $\rho_{\rm sep}=\mathcal E_{\rm sep}(f)/\norm{P_{[0,J]}f}$ is the error of \cref{thm:sep} per unit of input energy in the active scale window; $\rho_2$ is the second-layer error $\big(\sum_p\norm{\Psi'_{2,p}(L_\tau f)-L_\tau\Psi'_{2,p}f}^2\big)^{1/2}$ per unit of energy of the first-layer outputs in the second-layer windows ($j_2<j_1$); $e_\Pi/e$ is the covariance error of the dilation-averaged channel $\Psi'_{1,\Pi}$, with $\Pi$ uniform on $\{1,\frac87,\frac43,\frac85\}$, divided by that of the translation-averaged channel $\Psi_1'$.}\label{tab:ext}
\end{table}

\paragraph{Findings.} (i) \emph{Scale-separated channels.} $\rho_{\rm sep}$ remains flat in $J$ within $2\%$ under white inputs, whereas under pink inputs it decreases slowly. These inputs therefore do not attain the factor $\sqrt{J+1}$ of \cref{thm:sep}(a), as no growth in $\sqrt{J+1}$ occurs. The separated error, in absolute terms, reaches up to $2.4$ times the summed-channel error: scales combine sub-additively inside the modulus (\cref{ssec:num-typ}) in the summed channel, but not in the separated channels. (ii) \emph{Second layer.} The second-layer error per unit of energy amounts to $1.4$--$1.9$ times $\rho_{\rm sep}$ and stays within the factor $2$ of \cref{cor:depth2}, with no growth as $J$ increases. As the number of scale pairs increases, the second-layer energy $\sum_p\norm{\Psi'_{2,p}f}^2$ rises; it amounts to $1$--$8\%$ of the first-layer energy $\sum_j\norm{\Psi'_{1,j}f}^2$, consistent with the energy decay of wavelet scattering. (iii) \emph{Dilation averaging.} Whole coarsest cells must again make up the grid resampled by $D_r^{-1}$; thus $r=N_{\rm s}/N'$ is allowed when $N'$ is a multiple of $2^{11}$. The resulting nodes are $\{1,\frac87,\frac43,\frac85\}$. We choose the uniform probability measure $\Pi$, since \cref{cor:dil} permits any such measure. For each individual node, the error lies within $25\%$ of the $r=1$ value, in agreement with \cref{cor:dil}. Relative to the translation-averaged error, the averaged error $e_\Pi$ is $60$--$78\%$ and remains flat in $J$. Thus dilation averaging lowers the error further, presumably owing to the lack of alignment among errors at different nodes. A constant shift $\tau\equiv0.5$ gives a measured discretization floor of $6\times10^{-5}$.

\section{Lanczos versus the full matrix}\label{supp:lz}
The worst-case values of \cref{tab:worst} are Lanczos estimates of the largest eigenvalue of $A^*A$, where $A$ is the linear map of \cref{ssec:num-worst}. Lanczos alone cannot exclude a larger eigenvalue whose eigenvector is almost orthogonal to the starting vector. We therefore assembled the full $1024\times1024$ matrix of $A^*A$ on the band-limited space and computed all its eigenvalues, for Haar with $J\le3$ and aligned sym4 with $J\le2$ (smooth warp, $\eps=0.1$). In every case the largest eigenvalue agrees with the Lanczos value to a relative $8\times10^{-6}$ (\cref{tab:lzfull}).

\begin{table}[h]\centering\small
\begin{tabular}{lccccccc}\toprule
gen. & $J$ & steps & Lanczos $\theta^{1/2}/\eps$ & full $\lambda_{\max}^{1/2}/\eps$ & rel.\ diff.\ of $\lambda$ & $r/\theta$ & asymmetry\\\midrule
Haar & 0 & 40 & 1.294109 & 1.294109 & $2.1\times10^{-9}$ & $2.7\times10^{-7}$ & $4\times10^{-16}$\\
Haar & 1 & 32 & 1.507956 & 1.507956 & $1.3\times10^{-10}$ & $5.4\times10^{-10}$ & $5\times10^{-16}$\\
Haar & 2 & 40 & 1.944959 & 1.944959 & $1.6\times10^{-7}$ & $8.7\times10^{-6}$ & $5\times10^{-16}$\\
Haar & 3 & 40 & 1.977411 & 1.977419 & $8.0\times10^{-6}$ & $7.9\times10^{-4}$ & $5\times10^{-16}$\\
sym4 & 0 & 40 & 2.073101 & 2.073101 & $1.6\times10^{-8}$ & $1.5\times10^{-7}$ & $5\times10^{-16}$\\
sym4 & 1 & 24 & 2.387560 & 2.387560 & $1.1\times10^{-10}$ & $4.8\times10^{-10}$ & $4\times10^{-16}$\\
sym4 & 2 & 29 & 2.385739 & 2.385739 & $1.1\times10^{-10}$ & $4.8\times10^{-10}$ & $5\times10^{-16}$\\\bottomrule
\end{tabular}
\caption{Lanczos versus the full $1024\times1024$ matrix of $A^*A$ (smooth warp, $\eps=0.1$, $2^{13}$ samples). ``rel.\ diff.'' is $\abs{\theta-\lambda_{\max}}/\lambda_{\max}$, $r/\theta$ the relative Lanczos residual, and ``asymmetry'' is $\max\abs{M-M^T}/\max\abs M$ for the assembled matrix $M$.}\label{tab:lzfull}
\end{table}


\begin{thebibliography}{99}\small
\bibitem[BM13]{BM13} J.~Bruna, S.~Mallat, Invariant scattering convolution networks, \emph{IEEE Trans.\ Pattern Anal.\ Mach.\ Intell.} 35 (2013), 1872--1886.
\bibitem[Cal65]{Cal65} A.~P.~Calder\'on, Commutators of singular integral operators, \emph{Proc.\ Nat.\ Acad.\ Sci.\ U.S.A.} 53 (1965), 1092--1099.
\bibitem[CCM14]{CCM14} X.~Chen, X.~Cheng, S.~Mallat, Unsupervised deep Haar scattering on graphs, \emph{Advances in Neural Information Processing Systems} 27 (2014), 1709--1717.
\bibitem[Chu26]{Chu26} D.~Chung, Random Dyadic Shift Scattering: code and data, Zenodo (2026), \url{https://doi.org/10.5281/zenodo.24518588}.
\bibitem[CD95]{CD95} R.~R.~Coifman, D.~L.~Donoho, Translation-invariant de-noising, in \emph{Wavelets and Statistics}, Lecture Notes in Statist.\ 103, Springer, 1995, 125--150.
\bibitem[CM75]{CM75} R.~R.~Coifman, Y.~Meyer, On commutators of singular integrals and bilinear singular integrals, \emph{Trans.\ Amer.\ Math.\ Soc.} 212 (1975), 315--331.
\bibitem[Dau92]{Dau92} I.~Daubechies, \emph{Ten Lectures on Wavelets}, SIAM, 1992.
\bibitem[DL92]{DL92} I.~Daubechies, J.~C.~Lagarias, Two-scale difference equations II. Local regularity, infinite products of matrices and fractals, \emph{SIAM J.\ Math.\ Anal.} 23 (1992), 1031--1079.
\bibitem[DJ84]{DJ84} G.~David, J.-L.~Journ\'e, A boundedness criterion for generalized Calder\'on--Zygmund operators, \emph{Ann.\ of Math.} 120 (1984), 371--397.
\bibitem[DP22]{DP22} K.~Domelevo, S.~Petermichl, The dyadic and the continuous Hilbert transforms with values in Banach spaces, preprint, arXiv:2212.00090 (2022).
\bibitem[DP26]{DP26} K.~Domelevo, S.~Petermichl, A two-sided linear estimate and a dyadic reduction of the UMD conjecture, preprint, arXiv:2604.11273 (2026).
\bibitem[DWW23]{DWW23} F.~Di~Plinio, B.~D.~Wick, T.~Williams, Wavelet representation of singular integral operators, \emph{Math.\ Ann.} 386 (2023), 1829--1889.
\bibitem[Gra14]{Gra14} L.~Grafakos, \emph{Modern Fourier Analysis}, 3rd ed., Grad.\ Texts in Math.\ 250, Springer, 2014.
\bibitem[H08]{H08} T.~Hyt\"onen, On Petermichl's dyadic shift and the Hilbert transform, \emph{C.\ R.\ Math.\ Acad.\ Sci.\ Paris} 346 (2008), 1133--1136.
\bibitem[H12]{H12} T.~Hyt\"onen, The sharp weighted bound for general Calder\'on--Zygmund operators, \emph{Ann.\ of Math.} 175 (2012), 1473--1506.
\bibitem[HL22]{HL22} T.~Hyt\"onen, S.~Lappas, The dyadic representation theorem using smooth wavelets with compact support, \emph{J.\ Fourier Anal.\ Appl.} 28 (2022), 66.
\bibitem[M12]{M12} S.~Mallat, Group invariant scattering, \emph{Comm.\ Pure Appl.\ Math.} 65 (2012), 1331--1398.
\bibitem[NTV03]{NTV03} F.~Nazarov, S.~Treil, A.~Volberg, The $Tb$-theorem on non-homogeneous spaces, \emph{Acta Math.} 190 (2003), 151--239.
\bibitem[NT23]{NT23} F.~Nicola, S.~I.~Trapasso, Stability of the scattering transform for deformations with minimal regularity, \emph{J.\ Math.\ Pures Appl.} 180 (2023), 122--150.
\bibitem[Oja98]{Oja98} H.~Ojanen, Orthonormal compactly supported wavelets with optimal Sobolev regularity, arXiv:math/9807089.
\bibitem[P00]{P00} S.~Petermichl, Dyadic shifts and a logarithmic estimate for Hankel operators with matrix symbol, \emph{C.\ R.\ Acad.\ Sci.\ Paris} 330 (2000), 455--460.
\bibitem[PTV02]{PTV02} S.~Petermichl, S.~Treil, A.~Volberg, Why the Riesz transforms are averages of the dyadic shifts?, \emph{Publ.\ Mat.} (2002), Vol.\ Extra, 209--228.
\bibitem[S01]{S01} I.~W.~Selesnick, Hilbert transform pairs of wavelet bases, \emph{IEEE Signal Process.\ Lett.} 8 (2001), 170--173.
\bibitem[SBK05]{SBK05} I.~W.~Selesnick, R.~G.~Baraniuk, N.~G.~Kingsbury, The dual-tree complex wavelet transform, \emph{IEEE Signal Process.\ Mag.} 22 (2005), no.~6, 123--151.
\bibitem[Vag10]{Vag10} A.~Vagharshakyan, Recovering singular integrals from Haar shifts, \emph{Proc.\ Amer.\ Math.\ Soc.} 138 (2010), 4303--4309.
\bibitem[WB18]{WB18} T.~Wiatowski, H.~B\"olcskei, A mathematical theory of deep convolutional neural networks for feature extraction, \emph{IEEE Trans.\ Inform.\ Theory} 64 (2018), 1845--1866.
\end{thebibliography}
\end{document}